\documentclass{article}
\usepackage{latexsym,amssymb,epsf,amsmath,amsfonts,theorem}

\title{Stability of Poisson Equilibria and Hamiltonian\\
Relative Equilibria by Energy Methods}
\date{January 2002}
\author{George W.\ Patrick\\[1ex]
Mathematics and Statistics\\
University of Saskatchewan\\
Saskatoon, Saskatchewan, S7N~5E6\\
Canada \\
patrick@math.usask.ca
\and
Mark Roberts\\[1ex]
Mathematics and Statistics\\
University of Surrey\\
Guildford GU2 7XH\\
United Kingdom\\
m.roberts@surrey.ac.uk
\and
Claudia Wulff\,\thanks{On leave from: Freie Universit\"at Berlin}\\[1ex]
Mathematics Institute\\
University of Warwick\\
Coventry CV4 7AL\\
United Kingdom\\
claudia@maths.warwick.ac.uk}

\newtheorem{theorem}{Theorem}[section]
\newtheorem{lemma}[theorem]{Lemma}
\newtheorem{proposition}[theorem]{Proposition}
\newtheorem{corollary}[theorem]{Corollary}
\newtheorem{remark1}[theorem]{Remark}
\newenvironment{remark}{\begin{remark1}\rm}{\end{remark1}}
\newtheorem{definition1}[theorem]{Definition}
\newenvironment{definition}{\begin{definition1}\rm}{\end{definition1}}
\newtheorem{example1}[theorem]{Example}
\newenvironment{example}{\begin{example1}\rm}{\end{example1}}

\newcommand{\ro}{{\rm o}}
\newcommand\mfk\mathfrak
\newcommand\mcl\mathcal
\newcommand\mbb\mathbb
\newcommand\mtl\mathit
\newcommand\mbf\mathbf
\newcommand\onm\operatorname
\newcommand{\set}[2]{\left\{\,#1:#2\,\right\}}
\newcommand{\SE}{\mtl{SE}}
\newcommand{\SO}{\mtl{SO}}
\newcommand{\se}{\mtl{se}}
\newcommand{\so}{\mtl{so}}
\newcommand{\spl}{\mtl{sl}}
\newcommand{\g}{\mfk{g}}
\newcommand{\n}{\mfk{n}}
\newcommand{\mfsl}{\mfk{sl}}
\newcommand{\R}{\mbb{R}}

\newcommand{\V}{\mbb{V}}
\newcommand{\ad}{\onm{ad}}
\newcommand{\Ad}{\onm{Ad}}

\newcommand{\ink}{\rule{.3\baselineskip}{.35\baselineskip}}
\newcommand{\bpr}{\begin{trivlist} \item[]{\bf Proof. }}
\newcommand{\epr}{~~\nolinebreak\hspace*{\fill}$\ink$\end{trivlist}}

\begin{document}\maketitle

\begin{abstract}
\noindent
We develop a general stability theory for equilibrium points of
Poisson dynamical systems and relative equilibria of Hamiltonian
systems with symmetries, including several generalisations of the
Energy-Casimir and Energy-Momentum methods. Using a topological
generalisation of Lyapunov's result that an extremal critical point of
a conserved quantity is stable, we show that a Poisson equilibrium is
stable if it is an isolated point in the intersection of a level set of a
conserved function with a subset of the phase space that is related to
the non-Hausdorff nature of the symplectic leaf space at that
point. This criterion is applied to generalise the Energy-Momentum
method to Hamiltonian systems which are invariant under non-compact
symmetry groups for which the coadjoint orbit space is not Hausdorff.
We also show that a $G$-stable relative equilibrium satisfies the
stronger condition of being $A$-stable, where $A$ is a specific
group-theoretically defined subset of $G$ which contains the momentum
isotropy subgroup of the relative equilibrium.
\end{abstract}

\pagebreak
\footnotesize
\tableofcontents
\normalsize

\section{Introduction}

Energy methods for determining the stability of trajectories
of Hamiltonian systems are based on the general
principle that an
equilibrium point is Lyapunov stable if it is a strict local minimum
or maximum of a conserved function, such as the Hamiltonian itself.
The \emph{Energy-Casimir} and \emph{Energy-Momentum methods} are
extensions of
this principle to, respectively, equilibrium points of Poisson systems
and relative equilibria of Hamiltonian systems with symmetry.
They originate with the work of Arnold on the stability of equilibria
of incompressible fluids \cite{Arnold69}. Since then
they have been used very extensively in applications to
rigid bodies \cite{LeonardNEMarsdenJE-1997.1, LRSM92, P89},
elasticity theory \cite{LS90, PosberghTASimoJCMarsdenJE-1989.1,
SimoJCPosberghTAMarsdenJE-1989.1,SPM91},
fluids \cite{ArnoldVI-1978.1,AK98,CM90,
HolmDDMarsdenJERatiuTSWeinsteinA1985,MP94},
and vortex structures \cite{LP02,P00,PM98,W86}.
In this paper we first present a topological generalisation
of the energy method
and then use this to obtain
significant generalisations of the Energy-Casimir and Energy-Momentum
methods.

The flow of a Poisson system on a Poisson manifold $X$ generated by a
Hamiltonian $h$ preserves both $h$ and the symplectic leaves of $X$.
A point $x_e$ is an equilibrium point
of this flow if and only if it is a critical point of the restriction
of $h$ to the leaf $L(x_e)$ through $x_e$. If $x_e$ is a local
extremum of the restriction then the standard energy method implies
that $x_e$ is stable as an equilibrium point of the flow on $L(x_e)$.
In this case we say that $x_e$ is \emph{leafwise stable.} In
general $x_e$ is not a critical point of $h$ on the full space $X$.
To test for stability on the whole of $X$ the
Energy-Casimir method supposes that there is a
function $C$, the Casimir, which is constant on symplectic leaves
and such that $x_e$ is a critical point of $h+C$.
Stability follows if this critical point is a local extremum.

When is it possible to find a Casimir $C$ such that $x_e$ is a critical point
of $h+C$? One case is when $x_e$ is a \emph{regular} point of $X$, which
means that locally the foliation into symplectic leaves is non-singular.
Using this fact Arnold \cite{ArnoldVI-1978.1} and
Libermann and Marle \cite{LibermannPMarleCM-1987.1}
show that if $x_e$ is
regular and is a local extremum of the restriction
of $h$ to the leaf $L(x_e)$, then $x_e$ is stable for the full flow on
$X$. Thus at regular points this test for leafwise stability is also a test
for full stability. Examples show that this is not true in general, see
\cite[Exercise IV 15.10]{LibermannPMarleCM-1987.1} and Examples
\ref{sl(2)}, \ref{se(2)+} and \ref{RsdR} of this paper. In such cases
it is natural to ask whether there exists a space between $L(x_e)$ and $X$
such that $x_e$ is stable if it is an extremal point of the restriction of
$h$ to the intermediate space. In this paper we show that there is. More
generally we
answer a challenge posed by Weinstein \cite{WeinsteinA-1984.1} when,
referring to the interaction
between Poisson structures and stability, he wrote: ``As yet there is no
general theory for this kind of analysis''.

Most of the results in this paper are based on a topological generalisation
of the energy method (Corollary \ref{TopStabLem2}) which generalises a
lemma of Montaldi
\cite{MontaldiJA-1997.1}. Corollary \ref{TopStabLem2}
is valid for a continuous flow on a
locally compact topological space $X$ which has conserved quantities
with values in another topological space. In the case of Poisson
systems the conserved quantities are the Hamiltonian $h$ and the
quotient map to the space of symplectic leaves.
An equilibrium $x_e$ is stable if the leafspace is Hausdorff at $L(x_e)$ and
$x_e$ is an
isolated point in the fibre of the restriction of $h$ to $L(x_e)$. Thus the
condition that $x_e$ be regular in the result of Arnold, Libermann and Marle
can be relaxed to the leafspace being Hausdorff at $L(x_e)$. If the
leafspace is not Hausdorff then $h$ must isolate $x_e$ in a larger
subset $T_2(x_e)$ which depends only on the topology of the leafspace
(Theorem \ref{PoiStabTh1}).

We recover and generalise the Energy-Casimir method for Poisson
equilibria:
we will see that it suffices to make the assumptions of the
Energy-Casimir method on a subset of the Poisson manifold $X$ which
contains $T_2(x_e)$.
Example~\ref{3planes} shows that this improvement can succeed where
the standard Energy-Casimir method fails.
Moreover we identify a necessary condition for the Energy-Casimir method to
apply, namely that the Poisson equilibrium must be {\em tame}.
Our topological results provide stability tests which are even more general.

We now turn to the stability of relative equilibria. We consider
a $G$-invariant Hamiltonian $H$ on a symplectic manifold $P$ with a
momentum map $J:P\to \g^*$ which commutes with the action of $G$ on
$P$ and its coadjoint action on $\g^*$, the dual of the Lie algebra
$\g$ of $G$.
If the
action of $G$ on $P$ is free and proper then the orbit space $P/G$ is a Poisson
manifold and criteria for the stability of Poisson equilibria can be
lifted to criteria for the $G$-stability of relative equilibria. The
symplectic leaves of $P/G$ are just the Marsden-Weinstein reduced
phase spaces and $J$ induces a
homeomorphism between the leafspace of $P/G$ and the
coadjoint orbit space $\g^*/G$.

Leafwise stability of a relative equilibrium $p_e\in P$ means that it is stable
to {\em momentum-preserving} perturbations of the initial condition,
and is implied if the relative equilibrium is an extremal point of the
reduced Hamiltonian.
If the momentum $\mu_e = J(p_e)$ is
regular, or more generally $\g^*/G$ is Hausdorff at $G\mu_e$, then this
condition also implies that $p_e$ is $G$-stable
(Remark \ref{leafstab}).
When the Hessian of the reduced Hamiltonian is definite and the
momentum regular this was proved
by Libermann and Marle \cite{LibermannPMarleCM-1987.1},
following earlier work of Arnold \cite{ArnoldVI-1978.1}
and Marsden and Weinstein \cite{MW74}. The result for $\g^*/G$ Hausdorff
is due to Montaldi \cite{MontaldiJA-1997.1}.

The Energy-Momentum method is a convenient lifting of these criteria to
the phase space. A point $p_e$ is a relative equilibrium if and only if it
is a critical
point of the energy-momentum function $H_{\xi_e} = H - J_{\xi_e}$,
where $\xi_e$ is the {\em generator} of the relative equilibrium, an
element of the Lie algebra of $G$, and $J_{\xi_e}(p) = \langle
J(p),\xi_e\rangle$. The relative equilibrium is leafwise stable if the
restriction of the Hessian $d^2H_{\xi_e}(p_e)$ to a {\em symplectic normal
space}
$N_1$ is definite
\cite{MarsdenJESimoJCLewisDPosberghTA-1989.1,PosberghTASimoJCMarsdenJE-1989.1,
SLM91,SimoJCPosberghTAMarsdenJE-1989.1,SPM91}.
Full stability results were obtained by
Patrick \cite{P92} for compact groups,
and Ortega and Ratiu \cite{OrRa99} and Lerman
and Singer \cite{LeSi98} under more general assumptions which still
imply that $\g^*/G$ is Hausdorff at $G\mu_e$. See also Corollary
\ref{HamHessCor}.

Some stability results for non-compact groups at momentum values $\mu_e$ where
$\g^*/G$ is not Hausdorff have been obtained by Leonard and Marsden
\cite{LeonardNEMarsdenJE-1997.1} for semidirect products of compact
groups and vector spaces.
They suggested that it is necessary to test for the definiteness of
$d^2H_{\xi_e}(p_e)$ on a larger subspace of $T_{p_e}P$ than $N_1$. In this
paper we
sharpen and generalise their results to arbitrary groups by
lifting the general stability criteria for Poisson equilibria to conditions
on $H_{\xi_e}$. In particular we identify subsets containing
$N_1$ on which the definiteness of $d^2H_{\xi_e}(p_e)$
is sufficient to imply $G$-stability (Theorem \ref{HamHessStabTh1}
and Corollary \ref{Cor:Proceedings}).

The papers \cite{LeSi98,OrRa99,P92} all show, under various conditions
which in particular
imply that $\g^*/G$ is Hausdorff at $G\mu$, that the definiteness of
$d^2H_{\xi_e}(p_e)$ on $N_1$ implies that $p_e$ is actually $G_{\mu_e}$-stable,
ie trajectories that start near $p_e$ remain near $G_{\mu_e} p_e$, where
$G_{\mu_e}$
is the coadjoint isotropy subgroup at $\mu_e$. However numerical integration
of an example by Leonard and Marsden \cite{LeonardNEMarsdenJE-1997.1}
suggests that this is
not always true, even for a regular momentum value $\mu_e$. They
prove that these relative equilibria are $\Gamma$-stable for a subgroup
$\Gamma$ which lies strictly between $G_{\mu_e}$ and $G$.

In this paper we introduce
the general notion of $A$-stability for any subset of $G$ and then show
that any $G$-stable relative equilibrium is automatically $A$-stable
where $A$ is a `cone about $G_{\mu_e}$' which can be made arbitrarily close to
$G_{\mu_e}$ by restricting the perturbations from $p_e$ to be sufficiently
small
(Theorem \ref{AStabThm}). In many cases this result can be improved by
decomposing $G_{\mu_e}$
into the product of an (essentially) compact subgroup and a non-compact
submanifold
and showing that the `cone' only needs to be taken about the non-compact part.
(Theorem \ref{splitmu_A}). As a corollary we obtain a generalisation of
the results of \cite{LeSi98,OrRa99,P92} on $G_{\mu_e}$-stability
(Corollary \ref{compact_A}). Unlike the previous results
we do not require a Hessian condition to be satisfied, only that the
relative equilibrium is $G$-stable.

Taken together our results on the stability of Poisson equilibria and the
$G$ and $A$-stability of relative equilibria provide
generalisations of all previous results for finite dimensional Poisson
and Hamiltonian systems that we are aware of.
We believe that our
most primitive topological tests for
stability (Theorems \ref{PoiStabTh1} and \ref{HamTopThm}) are the sharpest
possible general
results, and that they fully explain the interaction between stability
properties and Poisson structures.
In this paper we
restrict attention to free Hamiltonian group actions for simplicity only;
the extensions to general proper actions are contained in \cite{PRW2}.

We end this introduction with a guide to the paper. 
Section \ref{TopologicalSect} contains the
fundamental topological ideas that underly the rest of the paper.
Topological and derivative tests for the stability of Poisson
equilibria are described in Sections
\ref{PoiTopSect} and \ref{PoiDerSect}, respectively.
Section \ref{PoiDerSect} introduces the machinery
of \emph{smoothings} that is used to make sense of derivatives of functions
on singular sets and the related notion of {\it tame generators}. It also
describes the role Casimirs play in stability theory and their limitations.
The topological and derivative tests for the $G$-stability of
Hamiltonian relative equilibria are contained in Sections
\ref{HamStabTestSect} and \ref{T2EMCM}, respectively.
The intermediate Section \ref{transverse} introduces
\emph{transverse Poisson structures} as tools for describing local
leafspace topology and uses this to discuss some special cases.
Towards the end of Section \ref{T2EMCM} we discuss how the stability criteria apply
to relative equilibria of systems that are invariant under actions of
the Euclidean groups $\SE(2)$ and $\SE(3)$ and how our theory
`explains' an example of 
Libermann and Marle \cite{LibermannPMarleCM-1987.1}.
Section \ref{AStability} is devoted
to $A$-stability: Section \ref{AStabilityThm} gives the general result,
Section \ref{Asplit}
the improvements obtained using Lie group decompositions and Section
\ref{AEuclidean} the application to Euclidean groups.

\section{Topology and Stability}
\label{TopologicalSect}

Let $X$ and $Y$ be a topological spaces and $f\colon X\rightarrow Y$ a
continuous map. In this section we give a stability criterion for
equilibria of continuous flows on X which preserve the fibres of $f$.
This is a corollary of a `confinement' result for general continuous
curves preserving $f$. The confinement result is inspired by Lemma 1.4
of Montaldi ~\cite{MontaldiJA-1997.1} in which $Y$ is Hausdorff.

We begin with
a measure of the extent to which a space is not Hausdorff at a point.
\begin{definition}\label{T2Df}
Let $Y$ be a topological space and $y\in Y$. Define
\begin{equation*}
T_2(y)\equiv\set{y^\prime\in Y}
{\ U \cap U^\prime \neq \emptyset \mbox{ for all neighbourhoods
 $y \in U \subseteq Y$ and $y^\prime \in U^\prime \subseteq Y$}}.
\end{equation*}
\end{definition}
We say that a topological space $X$ is {\em Hausdorff at $x\in X$}
if $T_2(x) = \{x\}$.
The following is frequently useful and easily proved.
\begin{proposition}\label{T2ProdLem}\mbox{}
\begin{enumerate}
\item Let $X$ and $Y$ be topological spaces and $(x,y)\in X\times Y$.
Then $T_2(x,y)=T_2(x)\times T_2(y)$.
\item If $X$, $Y_1$ and $Y_2$ are topological spaces and 
$f_i\colon X\rightarrow Y_i$, $i=1,2$, are continuous, then $(f_1\times
f_2)^{-1}\bigl(T_2(y_1,y_2)\bigr)=f_1^{-1}\bigl(T_2(y_1)\bigr)\cap
f_2^{-1}\bigl(T_2(y_2)\bigr)$. In particular, if $Y_2$ is Hausdorff then
$(f_1\times f_2)^{-1}\bigl(T_2(y_1,y_2)\bigr)$ is the $y_2$ level set
of $f_2|f_1^{-1}\bigl(T_2(y_1)\bigr)$.
\end{enumerate}
\end{proposition}
The next result gives a sufficient condition for continuous
curves in $X$ preserving $f$ and starting near a point $x$
to be confined to small neighbourhoods of $x$.

\begin{lemma}{\rm(Topological Stability Lemma)}\label{TopStabLem1}
Let $X$ and $Y$ be topological spaces, $f\colon X\rightarrow Y$
a continuous map, $x\in X$ and $y= f(x)$. Assume that:
\begin{enumerate}
\item $X$ is locally compact at $x$;
\item There exists a neighbourhood $U$ of $x$ in $X$ such that
$f^{-1}\bigl(T_2(y)\bigr) \cap U = \{x\}$.
\end{enumerate}
Then for every neighbourhood $U$ of
$x$ there is a neighbourhood $V$ of $x$ such that if
$c\colon [0,1]\to X$ is a continuous curve for which $f\circ c$ is
constant and $c(0)\in V$, then $c(t)\in U$ for all $t \in [0,1]$.
\end{lemma}
If the second condition holds we say that
{\em $x$ is an isolated point of
 $f^{-1}\bigl(T_2(y)\bigr)$}.

\bpr
We will prove the lemma by showing that if $f$-constant curves are
not confined to small neighbourhoods of $x$ then $x$ is not an isolated
point in $f^{-1}\bigl(T_2(y)\bigr)$.
Let $\mcl B$ be a neighbourhood base at $x$ consisting of compact
sets.
For any $U\in\mcl B$ and any $V\in\mcl B$ such that $V\subseteq U$,
let $c_{U,V}\colon [0,1]\rightarrow X$ be a continuous curve
such that $f\circ c$ is constant, $c_{U,V}(0)\in V$, and
$c_{U,V}(1)\not\in U$. By connectedness of $[0,1]$
there exists $t_{U,V}\in(0,1]$ such that
$x_{U,V} = c_{U,V}(t_{U,V})\in\partial U$.

For fixed $U\in\mcl B$ the set $\{x_{U,V}\}_{V \in \mcl B}$ is a net
in $\partial U$ by reverse inclusion of the $V$'s.
Since $\partial U$ is compact there is a subnet, $\{x_{U,V_\lambda}\}$,
which converges to a point $z_U\in\partial U$. The continuity
of $f$ implies that $\{f(x_{U,V_\lambda})\}$ converges to $f(z_U)$.
Since the curves $c_{U,V}$ preserve $f$ we also have:
\begin{equation*}
f(x_{U,V_\lambda})
=f\bigl(c_{U,V_\lambda}(t_{U,V_\lambda})\bigr)
=f\bigl(c_{U,V_\lambda}(0)\bigr).
\end{equation*}
The net $c_{U,V_\lambda}(0)$ converges to $x$ and so $f(x_{U,V_\lambda})$
also converges to $y=f(x)$. Thus, every neighbourhood of $f(z_U)$
meets every neighbourhood of $y$, and so $f(z_U)\in T_2(y)$.
This is the required contradiction, since
$z_U\in f^{-1}\bigl(T_2(y)\bigr)$ and $\{z_U\}_{U \in \mcl B}$ is a
net converging to $x$.
\epr

Now consider a continuous flow $\phi\colon X\times\mbb R\rightarrow X$
on $X$ which preserves the fibres of $f$,
$f\bigl(\phi_t(x)\bigr) = f(x)$,
and which has an equilibrium point at $x\in X$. The equilibrium
point is \emph{stable} if for every open neighbourhood $U$ of
$x$ in $X$ there exists a neighbourhood $V$ of $x$
such that if $x\in V$ then
$\phi_t(x) \in U$ for all $t$.

\begin{corollary}\label{TopStabLem2}
Let $X$ and $Y$ be topological spaces, $X$ locally compact,
and $f\colon X\rightarrow Y$ a continuous map. Let $x$ be an equilibrium
point of a continuous flow on $X$ which preserves $f$ and $y = f(x)$.
Then $x$ is stable if it is an isolated point in $f^{-1}\bigl(T_2(y)\bigr)$.
\end{corollary}

\noindent
The examples in Section~\ref{examples} show that Corollary~\ref{TopStabLem2}
is false in general if $x$ is only isolated in $f^{-1}(y)$.
The following example shows that it is necessary to assume that $X$ is
locally compact:

\begin{example}\label{20}
For $x=(q_n,p_n)_{n\ge0}\in X=\ell^2(\R)\times\ell^2(\R)$ let
\begin{equation*}
h(x)=\frac12\sum_{n=0}^\infty p_n^2
+\frac12\sum_{n=0}^\infty 4^{-n}q_n^2.
\end{equation*}
Then $h$ is differentiable on $X$ and $h(0)=0$. Clearly $T_2(0)=\{0\}$
and $0$ is isolated in $h^{-1}(0)\subset X$. The standard linear
symplectic structure $J\bigl((q_n,p_n)\bigr)=(p_n,-q_n)$ gives the
Hamiltonian system $\dot{x} = J d h(x)$ and defines a flow $\phi_t$ on
$X$. For $m\geq1$ the solution of this Hamiltonian system with
initial value such that $q_m=2^{-m}$, $q_n=0$ for $n\ne m$ and $p_n=0$
for all $n$, has $p_m=-1$ at time $t=2^{m-1}\pi$, as is easily verified.
Hence $0$ is an unstable equilibrium.
\end{example}

The requirement that $X$ is locally compact implies that
Corollary~\ref{TopStabLem2} can not be applied directly to deduce
the nonlinear stability of equilibria of partial differential equations.
Indeed, it is well-known that for partial differential equations
positivity of the second variation need not
imply stability of an equilibrium, a phenomenon which is related to
the non-equivalence of norms on Banach spaces \cite{BallMarsden}.

\section{Stability of Poisson Equilibria}
\label{PoiStabSect}

In this section we apply the topological stability result of Section
\ref{TopologicalSect} to equilibria of Poisson systems to
obtain generalisations of the Energy-Casimir method.

\subsection{Topological Tests}\label{PoiTopSect}

A finite dimensional Poisson manifold $X$ is
partitioned into immersed submanifolds by its symplectic leaves, which
without loss of generality we will always assume to be connected.
Define two points in $X$ to be equivalent if they belong to the same
symplectic leaf, let $Z$ be the quotient of $X$ by this equivalence
relation and $L:X \to Z$ the quotient map. We will regard the
symplectic leaf $L(x)$ through $x$ as both a subset of $X$ and as a
point in $Z$.

A function $h:X\rightarrow\mbb R$ generates a vector field on $X$
which is uniquely defined by the requirement that $\dot f = \{h,f\}$
for all differentiable functions $f$ on $X$. The flow~$\phi_t$ of $X$
preserves the fibres of both $h$ and $L$, and hence those of the
product map $f = L\times h:X\to Y=Z\times\mbb R$. We can therefore
apply the topological stability result Corollary \ref{TopStabLem2}.

First note that $x_e\in X$ is an equilibrium point of $\phi_t$ if and
only if it is an equilibrium point of the restriction of the flow to
the invariant submanifold $L(x_e)$ or, equivalently, if and only if
the restriction of $h$ to $L(x_e)$ has a critical point at $x_e$. An
equilibrium $x_e$ is said to be \emph{leafwise stable} if it is stable
for the restricted flow on $L(x_e)$. Clearly stability in the full
space $X$ implies leafwise stability.

For $x \in X$, let $T_2(x) = L^{-1}\bigl(T_2(L(x))\bigr)$. Note
that $L(x) \subseteq T_2(x)$. More generally, for any open
neighbourhood $ U$ of $x$ in $X$ let $L^U(x)$ denote the symplectic
leaf of $U$ through $x$. This is the connected component of $L(x)\cap
U$ which contains $x$. Denote the space of these symplectic leaves by
$Z^U$ and the corresponding quotient map by $L^U:U \to Z^U$. Define
\begin{equation*}
T_2^U(x)=(L^U)^{-1}\left(T_2(L^U(x))\right),
\end{equation*}
where $T_2\bigl(L^U(x)\bigr)$ is taken in $Z^U$. The set $T^U_2(x)$ is
contained in $T_2(x)\cap U$, but may be strictly smaller, for example
if $L(x)$ accumulates on itself at $x$.

Recall that a point $x_e$ in a Poisson manifold $X$ is said
to be {\em regular} (or {\em minimal}) if there exists an open
neighbourhood $U$ of $x_e$ in $X$ such that $\dim L(x)=\dim L(x_e)$
for all $x$ in $U$. The set of regular points is open and dense in
$X$ since in local coordinates it corresponds to the set where the
matrix of the Poisson tensor has locally constant rank.

The following result provides topological conditions for leafwise stability
and stability and states that for regular equilibria these conditions coincide.

\begin{theorem}\label{PoiStabTh1}
Let $x_e$ be an equilibrium point of the flow generated by a Hamiltonian
$h$ on a Poisson manifold $X$, and $U$ be an open neighbourhood of $x_e$.
Then $x_e$ is
\begin{enumerate}
\item
\emph{leafwise stable} if there is an open neighbourhood
$U$ of $x_e$ in $X$ such that
$h^{-1}\bigl(h(x_e)\bigr)\cap L^U(x_e)=\{x_e\}$; and 
\item
\emph{stable} if there is an open neighbourhood
$U$ of $x_e$ in $X$ such that
$h^{-1}\bigl(h(x_e)\bigr)\cap T^U_2(x_e)=\{x_e\}$.
\end{enumerate}
Moreover, if $x_e$ is a regular point of $X$ then there is an
open neighbourhood $U$ of $x_e$ in $X$ such that $T_2^U(x_e)=L^U(x_e)$
and part 1 implies stability.

These statements remain true if $h$ is replaced by any conserved quantity
with values
in a Hausdorff space.
\end{theorem}
\bpr

\noindent{\em 1.}\quad
This follows from Corollary \ref{TopStabLem2} with
$X$ replaced by $L^U(x_e)$ and $f$ by $h$. Clearly,
$T_2\bigl(h(x_e)\bigr)=\{h(x_e)\}$ and the hypothesis says that $x_e$ is an
isolated point in the $h(x_e)$ level set of the restriction of $h$ to $L^U(x_e)$.

\smallskip\noindent{\em 2.}\quad
Apply Corollary \ref{TopStabLem2} to the
neighbourhood $U$ and map $f=L^U\times h$, noting that by
Proposition~\ref{T2ProdLem} we have
$T_2(y)=T_2\bigl(L^U(x_e),h(x_e)\bigr)
=T_2\bigl(L^U(x_e)\bigr)\times\{h(x_e)\}$ and so
$f^{-1}\bigl(T_2(y)\bigr)=h^{-1}\bigl(h(x_e)\bigr)
\cap T^U_2(x_e)$.

\medskip\noindent
If $x_e$ is regular then there exists an open neighbourhood
$U$ of $x_e$ in $X$ for which the symplectic leaves provide a regular
foliation \cite[Corollary 2.3]{WeinsteinA-1983.1}.
The quotient of $U$ by this foliation is Hausdorff and so
$T^U_2(x_e)$ is equal to $L^U(x_e)$.
\epr

In Example~\ref{sl(2)} we show that there is an open subset
of Hamiltonians on $\mfsl(2)^* \cong \R^3$ for which
Theorem~\ref{PoiStabTh1} implies that the origin is a stable
equilibrium, and another open set for which the origin is leafwise
stable but not stable.

The next result is a simple corollary of Theorem~\ref{PoiStabTh1} for a
very special case.

\begin{corollary}
\label{PoiStabTh2}
Let $x_e$ be an equilibrium point of a Hamiltonian $h$ on the Poisson
manifold $X$. Suppose that $T^U_2(x_e)$ is a one-dimensional submanifold
of $X$. Then $x_e$ is stable if $dh(x_e)$ is nonzero on
$T_{x_e}\bigl(T^U_2(x_e)\bigr)$.
\end{corollary}

\bpr
If $dh(x_e)$ is nonzero on $T_{x_e}\bigl(T^U_2(x_e)\bigr)$ then
the level set of $h$ containing $x_e$ intersects $T_2(x_e)$
transversely and hence in an isolated point.
\epr

See Example~\ref{se(2)} for an illustration of this result. The
hypotheses of Corollary~\ref{PoiStabTh2} can only be satisfied if $x_e$
is a zero-dimensional symplectic leaf, in which case $x_e$ is
trivially a leafwise stable equilibrium for any Hamiltonian. In
Example~\ref{se(2)+} we show that the assumption that $T^U_2(x_e)$ be
one-dimensional cannot be replaced by the assumption that $L(x_e)$ has
codimension one in $T^U_2(x_e)$. Results analogous to
Corollary~\ref{PoiStabTh2}, in the sense that $dh(x_e)\ne0$ and
stability can be determined through a condition on $dh(x_e)$, can also
be obtained for some cases in which $T_2(0)$ is higher dimensional but
singular at $x_e$. See Example~\ref{sl(2)}.

\subsection[$T_2$-Energy-Casimir Method]{$\mathbf{T_2}$-Energy-Casimir Method}
\label{PoiDerSect}

We now discuss conditions on the derivative and Hessian of $h$
at $x_e$ for stability to hold.
If $x_e$ is an equilibrium, and so the restriction of $h$ to
$L(x_e)$ has a critical point at $x_e$, then it is leafwise stable if
the second derivative of the restriction is definite. The following
result, a special case of Theorem~\ref{PoiStabTh1}, states that for
generic points in $X$ this condition also implies that $x_e$ is stable.

\begin{proposition}\label{PoiDerTh1}
{\rm \cite[Theorem~III.12.4]{LibermannPMarleCM-1987.1}}
If an equilibrium $x_e$ is regular and the second derivative of the
restriction of the Hamiltonian $h$ to $L(x_e)$ is positive or negative
definite, then $x_e$ is both leafwise stable and stable.
\end{proposition}

If $T^U_2(x_e)$ contains a two dimensional manifold passing
through $x_e$ on which $dh(x_e)\ne 0$ then the level set of $h$
in $T^U_2(x_e)$ through $x_e$ contains a smooth curve
and so $h$ does not isolate $x_e$ in $T^U_2(x_e)$.
In such a case stability cannot be concluded from the
mechanism of confinement by energy level sets. However, if
$T^U_2(x_e)$ is a manifold, if the first derivative of $h$ vanishes
on $T^U_2(x_e)$ and if the second derivative of the restriction of $h$
to $T^U_2(x_e)$ is definite, then $h$ again isolates $x_e$ and $x_e$
is stable.
We call Poisson equilibria for which the first derivative of $h$ satisfies
this condition, and which are therefore amenable to Hessian type stability
tests, {\em tame}.
To extend this to cases for which $T^U_2(x_e)$ is not a manifold,
we first define the notions of {\em tangent space} and \emph{smoothing}
of a singular set.

\begin{definition}\label{SmoothDf}
Let $M$ be a manifold, $S\subseteq M$ and $m\in S$.
\begin{enumerate}
\item
The \emph{tangent space} $T_mS$ of $S$ at $m$ is the subset of
$T_mM$ consisting of the derivatives $c^\prime(0)$ of all $C^1$
curves $c(t)$ in $M$ with $c(0)=m$ and $c(t) \in S$ for $t\ge0$.
\item
A \emph{smoothing} of $S$ at $m$ is a finite number of submanifolds
$B_i\subseteq M$ such that $S\subseteq\bigcup_{i=1}^n B_i$ and
$T_mB_i\subseteq\onm{span}(T_mS)$ for $1\le i\le n$.
\end{enumerate}
\end{definition}
If $T^U_2(x_e)$ is a manifold then clearly it is its own smoothing.
Examples \ref{sl(2)} and \ref{2planes} feature some smoothings of
singular $T^U_2(x_e)$ sets.
{\em Weak smoothings} which satisfy the first condition of 2., but not the
second,
always exist since it is possible to take $n=1$ and $B_1 = M$.
 In Example \ref{nosmoothing} the $T_2$-set does not have a smoothing.

We now define tame equilibria at points with arbitrary $T_2$-sets.
\begin{definition}\label{Ptame}\mbox{}
Let $X$ be a Poisson manifold.
\begin{enumerate}
\item
A \emph{generator} at $x$ is an element $\xi\in T^*_xX$ which
annihilates $T_xL(x)$.
\item
A generator~$\xi$ at $x$ is \emph{tame} if it annihilates
$T_{x}\bigl(T^U_2(x)\bigr)$. Generators that are not tame are said to
be \emph{wild}.
\item
The generator of an equilibrium $x_e$ of a Poisson system with
Hamiltonian $h$ is $dh(x_e)$.
\item
An equilibrium is tame if its generator is tame, and wild otherwise.
\end{enumerate}
\end{definition}
We will notationally suppress reference to $U$ even though the property of being
tame is $U$~dependent. The set of tame generators at $x$ is a vector
subspace of $T_x^*X$. If $\xi \in T_x^*X$ is tame then it annihilates
every smoothing $\left\{B_i\right\}_{i=1}^n $ of $T_2(x)$, ie $\xi$
annihilates $T_xB_i$ for all $i$.

If $x$ is regular then $T^U_2(x) = L^U(x)$ for some neighbourhood $U$
and every generator is tame. In Corollary~\ref{PoiStabTh2} the
hypothesis that $dh(x_e)$ is nonzero on $T_{x_e}\bigl(T^U_2(x_e)\bigr)$
implies that $dh(x_e)$ is wild.
In Example~\ref{sl(2)}, for which
$T_2(0)$ is a cone, $T_0\bigl(T_2(0)\bigr)$ spans the whole of $X
=\mbb R^3$ and so every nonzero derivative at $0$ is wild.

A~\emph{Casimir} on $U$ is a continuous function
$C:U\rightarrow\mbb R$ which is constant on the symplectic leaves of $U$.
This condition implies that $C$ is also constant on every set $T^U_2(x)$
for every $x\in U$. Casimirs are conserved quantities along integral curves
contained in $U$, since any Poisson flow preserves the
symplectic leaves. Since Casimirs are constant on $T^U_2(x)$,
the derivatives of smooth Casimirs at $x$ are tame.
It follows that if $x_e$ is an equilibrium and $C$ is a smooth Casimir
at $x_e$, then $dh(x_e)$ is tame if and only if $d(h+C)(x_e)$ is tame.

\begin{theorem}{\rm(Poisson $T_2$-Energy-Casimir Method)}
\label{PoiHessStabTh} Let $x_e$ be a tame equilibrium point of the
Hamiltonian $h$. Let $\{B_i\}_{i=1}^n$ be a smoothing of $T^U_2(x_e)$
at $x_e$. Then $x_e$ is
stable if for each $i$ there is a smooth Casimir $C_i$ such that the
Hessian $d^2\bigl((h+C_i)|B_i\bigr)(x_e)$ is positive or negative
definite on $T_{x_e}B_i$.
\end{theorem}

\bpr Set $\hat h = (h+C_1,\cdots,h+C_n)$.
Since $T^U_2(x_e)\subseteq\bigcup_{i=1}^n B_i$,
by Theorem~\ref{PoiStabTh1} it suffices to show that
$\hat h$ isolates $x_e$ on
$\bigcup_{i=1}^n B_i$.
As there are only finitely many $B_i$,
this follows if $h+C_i$ isolates $x_e$ on each $B_i$.
This in turn is implied by the definiteness of the Hessians and
the Morse lemma.\epr

The standard Energy-Casimir method, see eg \cite{MarsdenJERatiuTS-1994.1},
states that if $x_e$ is an equilibrium
point of the flow generated by a Hamiltonian $h$ and there exists a
smooth Casimir $C$ such that $x_e$ is a critical point of $h+C$ on the
whole of $X$ and $d^2(h+C)(x_e)$ is definite, then $x_e$ is stable.
Theorem~\ref{PoiHessStabTh} is a strict generalisation of this
because it only requires $x_e$ to be a definite critical point for a
function on a subset of $X$. Example~\ref{3planes} describes a Poisson
system for which the theorem can be used to prove stability, though
the standard Energy-Casimir method fails.

The standard Energy-Casimir method and its generalisation, Theorem
\ref{PoiHessStabTh}, can only be applied to tame equilibria,
but the topological stability result Theorem
\ref{PoiStabTh1} can sometimes be applied to wild equilibria.
Examples include the cases covered by Corollary
\ref{PoiStabTh2} and Examples~\ref{se(2)} and~\ref{sl(2)}.

Smoothings and Casimirs are both implements designed to handle
singularities of $T_2^U(x_e)$ for
the purpose of constructing Hessian tests for stability.
Fine smoothings can eliminate the
necessity of including Casimirs when calculating the Hessians, while
coarse smoothings may require their inclusion. Example~\ref{2planes}
illustrates the play available in choosing the smoothings versus
Casimirs.

If there exists a smoothing such that
$T_2^U(x_e)=\bigcup_{i=1}^{n}B_i$ (for example, if $T_2^U(x_e)$ is a
submanifold) then every Casimir is constant on each of the $B_i$. In
this case the second derivatives of the restrictions of the Casimirs
to the~$B_i$ vanish, and the inclusion of the Casimirs in the Hessians
in Theorem~\ref{PoiHessStabTh} is unnecessary. However,
Example~\ref{sl(2)} shows that inclusion of Casimirs can be necessary
in cases when $T^U_2(x_e)$ is singular.

More generally, the inclusion of the Casimir $C_i$ is unnecessary when
all Casimirs restricted to $B_i$ have vanishing second derivative at
$x_e$. This follows if $B_i\subseteq T_2^U(x_e)$ but can also be
implied by an infinitesimal relation between the smoothing
and the set $T_{x_e}^U(x_e)$, as follows. For any Casimir $C$ we have
$d^2C(x_e)(v,v)=0$ for all $v\in T_{x_e}\bigl(T_2^U(x_e)\bigr)$, since
Casimirs are constant on curves in $T_2^U(x_e)$. Regarding the
symmetric bilinear form $d^2C(x_e)$ as a linear map on the tensor
product $T_{x_e}B_i\otimes T_{x_e}B_i$, if follows that $d^2C(x_e)$
vanishes on the whole of $T_{x_e}B_i$ if
\begin{equation}\label{casUnnecEq}
\onm{span}\set{v\otimes v}{v\in T_{x_e}\bigl(T_2^U(x_e)\bigr)
\cap T_{x_e}B_i}=T_{x_e}B_i\otimes T_{x_e}B_i.
\end{equation}
Consequently, in Theorem~\ref{PoiHessStabTh}, the inclusion of $C_i$ is
unnecessary for any $i$ such that this spanning condition holds. An
application of this is given in Example~\ref{3planes}.

\begin{remark}\label{weakSmoothRem}
If $\left\{B_i\right\}_{i=1}^n$ is only a weak smoothing of
$T_2(x_e)$, but there exist Casimirs $C_i$
such that for $i=1,\ldots,n$ we have $d (h+ C_i)(x_e)|_{B_i}=0$ and
the Hessians $d^2 (h+ C_i)(x_e)|_{B_i}$ are definite, then $x_e$ is
again stable. This can be deduced as in
the proof of Theorem \ref{PoiHessStabTh}.
The standard Energy-Casimir
method is recovered by taking $n=1$ and $B_1 = M$.
\end{remark}

\subsection{Examples}
\label{examples}

Here we collect together a number of examples to illustrate the
stability theory described above. Many of the examples are of
equilibria on Poisson spaces which are duals of Lie algebras,
for which the symplectic leaves are coadjoint orbits.
Other examples are Poisson structures on
$X=\R^3$ with Poisson brackets of the form:
\begin{equation}\label{Abracket}
\{h,f\}=\nabla A\cdot(\nabla h\times\nabla f)
\end{equation}
where $A = A(x,y,z)$ is a smooth function. The vector field generated
by a Hamiltonian $h$ is
\[
\dot{x}=\nabla A \times \nabla h.
\]
We may assume that $A(0) =
0$. For these structures $A$ is a Casimir since $\{ A, f\} =0$ for
all functions $h$, and the two dimensional symplectic leaves of $X$
are the connected components of the level sets of regular values of
$A$. Each critical point $x$ of $A$ is also a symplectic leaf since
there $\{h,f\}(x)=0$ for all smooth functions $h$ and $f$. For any
open neighbourhood $U$ of $0$ in $X$ the set $T_2^U(x)$ is
contained in the connected component of $A^{-1}(0)\cap U$ which contains
$0$. If $A = \frac{1}{2}(x^2 + y^2 + \epsilon z^2)$ then $X$ is
isomorphic to the dual of the Lie algebra, $\g^*$, of $\SO(3)$
($\epsilon = 1$), $\SE(2)$ ($\epsilon = 0$) or $SL(2,\R)$ ($\epsilon =
-1$) with their standard Poisson structures. These three
Poisson structures are described in~\cite{MarsdenJERatiuTS-1994.1}
and are also used in
~\cite{WeinsteinA-1984.1} to illustrate the interdependence of
stability and Poisson structure.

\begin{example}
\label{so(3)}
If $A = \frac{1}{2}(x^2 + y^2 + z^2)$, and so $X \cong \so(3)^*$, then
for any neighbourhood $U$ of $0$ we have $L(0) = T_2^U(0) = \{0\}$ and
so $0$ is an isolated point of $h^{-1}\bigl(h(0)\bigr) \cap T_2^U(0)$, and hence
a stable equilibrium, for any Hamiltonian $h$. Every generator at $0$ is
tame. Every point $x_e \neq 0$ is regular, since the symplectic leaf
through each nearby point (a sphere) is two dimensional. Thus all
generators at these points are also tame. In fact every generator is
tame at every point in $\g^*$ for any compact group, since in this
case the symplectic leaves are the coadjoint orbits of $G$ on $\g^*$
and the quotient space $\g^*/G$ is always Hausdorff.
\end{example}

\begin{example}\label{se(2)}
If $A = \frac{1}{2}(x^2 + y^2)$ then $X \cong \se(2)^*$. The symplectic
leaves are points on the $z$-axis and cylinders of nonzero radius about
the $z$-axis.

If $x_e$ is a point on the $z$-axis then $T_2^U(x_e)$ is the
intersection of the $z$-axis with $U$, a one dimensional manifold.
The set of tame generators at $x_e$ is equal to the subspace
$\R^2 \subseteq \se(2)$ consisting of infinitesimal translations.
If $dh(x_e)$ is tame the point $x_e$ is stable if the restriction of
$h$ to the $z$-axis is positive or negative definite.
Generators containing nonzero rotational components are wild.
However in this case this implies that $dh(x_e)$ is nonzero on
$T_{x_e}\bigl(T_2^U(x_e)\bigr) = \{x=y=0\}$
and so $0$ is again stable by Corollary \ref{PoiStabTh2}.

If $x_e$ does not lie on the $z$-axis then it is regular and every generator
is tame.
\end{example}

\begin{example}\label{sl(2)}
Next consider $A = \frac{1}{2}(x^2 + y^2 - z^2)$,
giving $X \cong \spl(2,\R)^*$. The symplectic leaves are the connected
components of the hyperboloids $A = a$ for $a \neq 0$, the connected
components on the complement of $0$ in the cone $A=0$, and the origin
itself. Every nonzero point is regular and all generators at these
points are tame.

At the origin, $T_2^U(0)$ is the intersection of the cone
$A=0$ with $U$. Since $T_0\left(T_2^U(0)\right)$ spans the whole
of $\R^3$ the only tame generator is $\xi=0$. Nevertheless Hamiltonians
with wild generators $dh(0)\ne0$ can again have stable equilibria at $0$.
Theorem \ref{PoiStabTh1} implies that $0$ is stable if $dh(0)$ `points into
the cone $A = 0$',
so that $\onm{ann}\bigl(dh(0)\bigr)$ intersects the cone only at the origin. If
$dh(0)$ points out of the cone then the intersection is (infinitesimally)
a pair of lines
and the equilibrium is unstable. The instability follows from the fact
that if $dh(0) = (\xi_1,\xi_2,\xi_3)$ the linearised equations of motion
at the origin have eigenvalues $0$ and $\pm\sqrt{\xi_1^2+\xi_2^2-\xi_3^2}$.
Since $0$ is trivially leafwise stable for any Hamiltonian this is
an example of an equilibrium that is leafwise stable but not stable.

Now consider the tame case $dh(0)=0$. For a smoothing we have
to take $B = U$, a full neighbourhood of the origin. Without including
Casimirs we can only use Theorem \ref{PoiHessStabTh} to conclude stability if
$d^2h(0)$ is definite. However by using the Casimir $A$ it can be
seen that, for example, the Hamiltonian $ax^2+by^2+cz^2$ has a stable
equilibrium whenever $c > -a$ and $c > -b$. Thus in general Casimirs
can not be dispensed with completely.
\end{example}

\begin{example}\label{se(2)+}
Let $X$ denote the Poisson space $\se(2)^* \times \R^2$, with coordinates
$x,y,z$ on $\se(2)^*$, as in Example \ref{se(2)}, and coordinates $p,q$ on
$\R^2$ satisfying $\{p,q\} = 0$. Then the symplectic leaf through $0$
is the two dimensional manifold $L(0) = \{x=y=z=0\}$ and $T_2^U(0)$
is the intersection of $U$ with the three dimensional manifold
$\{x=y=0\}$. The set of generators can be identified with
$T^*_0(\se(2)^*) \cong \se(2)$ and the tame generators are the pure
translations, as in Example \ref{se(2)}.

However in contrast to Example~\ref{se(2)}, Corollary \ref{PoiStabTh2}
does not apply, and $0$ can be unstable when $dh(0)$ is wild. To see
this consider the Hamiltonian $h = az - qy + \frac{1}{2}(q^2 + p^2)$.
This has a wild generator at $0$. Its restriction to the symplectic leaf
has a strict local minimum at $0$, so $0$ is leafwise stable. The equations
of motion are:
\begin{equation*}
\dot x = ay, \quad
\dot y = -ax, \quad
\dot z = -qx, \quad
\dot q = p, \quad
\dot p = -q+y.
\end{equation*}
When $a=1$ for arbitrarily small $\delta$ the curve
\begin{equation*}
x=\delta\sin t,~~
y=\delta\cos t,~~
z=\frac{\delta^2}{8}(t\sin2t-\sin^2 t-t^2),~~
q=\frac\delta2t\sin t,~~
p=\frac\delta2(\sin t+t\cos t)
\end{equation*}
is a solution which leaves any neighbourhood of $0$ for $t$
sufficiently large. This instability is already present in the linearised
equations at $0$ and is caused by the $1:1$ resonance between the $x,y$
and $q,p$ frequencies.
\end{example}

\begin{example}\label{2planes}
Let $X = \R^3$ with the Poisson structure given by (\ref{Abracket}) with
$A = \frac{1}{2}(x^2 - y^2)$. The set $T_2^U(0)$ is the intersection
of $U$ with the union of the two planes $x = \pm y$.
The only tame generator is $\xi = 0$.
One smoothing of $T_2^U(0)$ is provided by the two planes themselves:
$B_1 = \{x=y\}$, $B_2 = \{x=-y\}$. Another is obtained by taking a whole
open neighbourhood $U$ of $0$.

Consider the Hamiltonian $h =ax^2 - by^2 + z^2$. The restriction of
$h$ to each $B_i$ is $(a-b)x^2+z^2$ and so $0$ is stable by
Theorem~\ref{PoiHessStabTh} if $a>b$, without using any Casimirs.
Alternatively the Casimir
$A$ can be used on the whole of $U$:
\begin{equation*}
h-2\lambda A=(a-\lambda)x^2-(b-\lambda)y^2+z^2
\end{equation*}
and so if $a>b$ then taking $a>\lambda>b$ gives stability by
the Energy-Casimir method.
\end{example}

\begin{example}\label{3planes}
This is an example of an equilibrium which is stable by Theorem
\ref{PoiHessStabTh}, but for which the standard Energy-Casimir method
fails. Again let $X = \R^3$ with the Poisson structure given by
(\ref{Abracket}), but with $A = (a^2x^2 - y^2)y$ where $a \neq 0$.
The set $T_2^U(0)$
is the intersection of $U$ with the union of the three planes $y = 0$
and $y = \pm ax$. The only tame generator is $\xi = 0$. The
restriction of the Hamiltonian $h = x^2 - y^2 + z^2$ to $y = 0$
isolates $0$ in $y =0$, while its restrictions to $y = \pm ax$ isolate
$0$ if $\vert a \vert \leq 1$. Thus $0$ is stable by Theorem
\ref{PoiHessStabTh} if $\vert a \vert \leq 1$. However any Casimir on
$X$ must satisfy $d^2C(0) = 0$ since, as $v$ varies over the three
planes forming the set $T_2(0)$, the vectors $v \otimes v$ span the
symmetric part of $\R^3 \otimes \R^3$. Thus the standard
Energy-Casimir method can not be used to deduce stability.
\end{example}

\begin{example}\label{Unecessary}
This example shows that the use of $T_2^U(x_e)$ instead of the larger
set $T_2(x_e)$ is necessary when the symplectic leaves accumulate upon
themselves. Let $\hat X =\R^3$ with the Poisson bracket~(\ref{Abracket})
with $A = x-ay$. Let
\begin{equation*}
\hat h=\frac12\bigl(z^2+(\sin y-a\sin x)^2-(\sin x-a\sin y)^2\bigr),
\end{equation*}
and let $\hat x_e = 0$. The action of $\mbb Z\times\mbb Z$ on $\hat X$ by
$(b_1,b_2)\cdot(x,y,z) = (x+2\pi b_1,y+2\pi b_2,z)$
is Poisson and $\hat h$ is invariant. Let $X\equiv\hat X/\mbb
Z\times\mbb Z$ be the Poisson quotient, and $h$ and $x_e$ be the
projections of $\hat h$ and $\hat x_e$ to the quotient. Let $a$ be
irrational. Then the symplectic leaf through $x_e$ is the projection of
the plane $x=ay$ in $\hat X$. This is dense in $X$ since it is the product
of a densely winding line on a 2-torus and $\mbb R$. So $T_2(x_e)=X$
and Theorem~\ref{PoiHessStabTh} fails to show that $x_e$ is stable,
since the Hessian of $h$ is indefinite on $X$.
However, by taking $U$ to be the projection of
$(-r,r)\times(-ar,ar)\times\mbb R$ for sufficiently small $r$,
the set $T^U_2(x_e)$
becomes the product of the projection of the \emph{line segment} $y=ax$,
$x\in(-r,r)$, and $\mbb R$, which is not dense. The
Hessian of $h$ restricted to $T^U_2(x_e)$ is definite since
$\hat h=\frac12(z^2+(1-a^2)^2y^2)+\mbox{\it h.o.t.}$ on the plane
$x=ay$ in $\hat X$.
\end{example}

\begin{example}\label{nosmoothing}
Let $X= \spl(3;\R)^*$. Using the Killing form we can identify $X$
with $\spl(3;\R)$, and hence with the space of traceless $3 \times 3$
real matrices. Let $x_e$ denote the \emph{subregular} nilpotent
matrix with
$(x_e)_{12} = 1$ and all other entries equal to zero.
We claim that $T_2^U(x_e)$ does not
have a smoothing for any neighbourhood $U$. The coadjoint orbit
$L(x_e)$ has codimension four. Let $\Sigma$ be a four-dimensional
section through $x_e$ transverse to $L(x_e)$.
This has an induced Poisson structure
that is described in Section \ref{transverse} below and
$T_2^U(x_e)$ is isomorphic to
the product of a neighbourhood of $x_e$ in $L(x_e)$ and
$T_2^{U_\Sigma}(x_e)$ where $U_\Sigma$ is a neighbourhood of $x_e$ in
$\Sigma$. We will show that $T_2^{U_\Sigma}(x_e)$ does not have a
smoothing for any neighbourhood $U_\Sigma$.

Let $\pi:X\to\R^2$ denote the mapping defined by generators of the
ring of invariants of the coadjoint action on $X$, ie the symmetric
polynomials of the eigenvalues of the matrices. Let $\pi_\Sigma$ be
the restriction of $\pi$ to $\Sigma$. Since the generators are
Casimirs the set $T_2^{U_\Sigma}(x_e)$ is contained in the
two-dimensional fibre of $\pi_\Sigma$ through $x_e$. Brieskorn has
shown that this fibre has a simple singularity of type $A_2$
\cite{B70,S80}, which means that it is diffeomorphic to the variety
defined by $x^2 + y^2 + z^3 = 0$. The tangent space to this variety
at $0$ (in the sense of Definition \ref{SmoothDf}) is just the
non-negative $z$-axis, ie
$T_{x_e}\bigl(T_2^{U_\Sigma}(x_e)\bigr)=\set{(0,0,z)}{z\ge 0}$. However
the fibre $\pi_\Sigma^{-1}(\pi_\Sigma(x_e))$ contains the intersection
of $\Sigma$ with the coadjoint orbit through the \emph{regular}
nilpotent matrix defined by $(x_e)_{12} = (x_e)_{23}= 1$ and all other
entries equal to zero. This intersection is a symplectic leaf of
$\Sigma$ of dimension two which, by Jordan normal form theory,
contains $x_e$ in its closure. It must therefore be contained in
$T_2^{U_\Sigma}(x_e)$. Thus a smoothing $\{B_i\}_{i=1}^{n}$ 
must have one $B_i$ at least of dimension~2 and it is not possible
that $T_{x_e}B_i\subset\mbox{span }T_{x_e}(T_2^{U_\Sigma}(x_e))$ 
for each $i$.
\end{example}

\begin{example}
\label{se(2)n}
This example shows that the assumption in Theorem~\ref{PoiHessStabTh} that
the equilibrium is tame is essential, even for Poisson systems that are
reductions of Hamiltonian systems of physically recognisable forms.
Take $P=T^*\bigl(\mtl{SE}(2)^n\bigr)$ and $G=\mtl{SE}(2)^n$ acting
by the cotangent lift of its left action on itself. Then $X = P/G$
is the Lie algebra dual $(\mtl{se}(2)^*)^n \cong (\mbb R^3)^n$,
and the generic coadjoint orbits are products of cylinders. In cylindrical
coordinates $(r_i,\theta^i,z_i)$ the Poisson bracket is
\begin{equation*}
\{f,g\}=\sum_{i=1}^n\left(
\frac{\partial f}{\partial z_i}\frac{\partial g}{\partial\theta^i}
 - \frac{\partial g}{\partial z_i}\frac{\partial f}{\partial\theta^i}
\right).
\end{equation*}
Consider Hamiltonians $h:P/G\to \R$ of the form
\begin{equation}\label{se2_blow}
h=\sum_{i=1}^n \frac12{z_i}^2+r_i^2 F_i(\theta ,z )
\end{equation}
for smooth functions $F_i$ on the product of $n$ cylinders,
$z = (z_1,\ldots, z_n)$, $\theta = (\theta^1,\ldots, \theta^n)$.
Then the corresponding Poisson system is given by
\begin{equation*}
\dot{r}_i=0,\quad\dot{z}_i = -\partial_{\theta^i} h,\quad\dot{\theta}^i =
\partial_{z_i}h,\quad i=1,\ldots, n.
\end{equation*}
For $i=1\ldots,n$ the coordinate functions $r_i$ are Casimirs
that parametrise the Hamiltonian~(\ref{se2_blow}).
Any point with all $r_i=0$ is a symplectic leaf, and hence a
(leafwise stable) equilibrium. The full stability of these equilibria
corresponds to stability under small perturbation to nonzero
$r_i$. These are perturbations of a completely integrable system,
and so when $n\ge3$ Arnold drift in the variables $z_i$ is
expected~\cite{ArnoldVI-1964.1}.
For $n=3$ Bessi~\cite{BessiU-1997.1} has constructed functions $F_i$
for which instability occurs near $z_1=z_3=0$ and $z_2=2$.
The $T_2$-set corresponding to any equilibrium with all $r_i=0$ is the set
$\set{(r_i,\theta^i,z_i)}{r_1=\cdots=r_n=0}$ and $dh$ restricted to
the $T_2$-set is not zero unless all $z_i$ are zero.
Thus the equilibrium is wild for $z \neq 0$.

We note that for $n=1$ stability follows from
Corollary~\ref{PoiStabTh2}, and that for $n=2$ we essentially have a
family of two-degree of freedom systems (the phase space is the product of two
cylinders of widths $r_1$ and $r_2$) so that generically stability follows
from KAM theory.
\end{example}

\section{Stability of Hamiltonian Relative Equilibria}
\label{HamStabSect}

In this section we apply the results of Section \ref{PoiStabSect}
on the stability of Poisson equilibria to relative equilibria
of Hamiltonian systems which are invariant under free group actions.
To do this we use the fact that relative equilibria are equilibria
of the symmetry reduced dynamics on the Poisson orbit manifold.

Let $P$ be a finite dimensional symplectic manifold with a free
symplectic action of a Lie group $G$.
Assume that the action has an $\Ad^*_G$-equivariant momentum map
$J:P\rightarrow\mfk g^*$ with respect to the coadjoint action of
$G$ on $\mfk g^*$.
Let $H$ be a $G$-invariant function on $P$,
$X_H$ the corresponding equivariant Hamiltonian vector field
and $\Phi_t$ its equivariant flow. The flow preserves the level sets of
both $H$ and $J$. Note that a fibre $J^{-1}(\mu)$ is only $G_\mu$-invariant
and the restriction of the flow to the fibre is $G_\mu$-equivariant.

By definition, a point $p_e\in P$ is a relative equilibrium if there
exists a {\it generator} $\xi_e\in\g$ such that $X_H(p_e)=\xi_ep_e$.
This is equivalent to $p_e$ being a critical point of
$H_{\xi_e}=H-J_{\xi_e}$ where $J_{\xi_e}(p)=J(p)(\xi_e)$.
Note that the trajectory of $X_H$ through $p_e$ is $\exp(t\xi_e)p_e$.
Equivariance and conservation of $J$ implies that $\xi_e \in \g_{\mu_e}$.
\begin{definition}\label{REStabDef}
A relative equilibrium $p_e$ is:
\begin{enumerate}
\item
{\it $G$-stable} if for every $G$-invariant
neighbourhood $U$ of $p_e$ there is a neighbourhood $V$ such
that $\Phi_t(p)\in U$ for all $p\in V$ and all $t$;
\item
{\it Leafwise stable} if it is $G_{\mu_e}$-stable for the
restriction of $\Phi_t$ to the momentum level set $J^{-1}(\mu_e)$.
\end{enumerate}
\end{definition}
Thus a relative equilibrium is leafwise stable if it is stable
(mod $G_\mu$) to momentum preserving perturbations and $G$-stable
if it is stable (mod $G$) to all perturbations.

\emph{For simplicity we restrict the discussion in this paper to free,
proper group actions.} In~\cite{PRW2} we will take advantage of the
generality of the topological stability lemma to extend the
stability theory to actions with nontrivial isotropy subgroups. For a
free, proper action, the orbit space $P/G$ is a smooth manifold
which inherits a Poisson structure from the symplectic structure on
$P$. The Hamiltonian $H$ descends to a function $h$ on $P/G$ for
which the corresponding Poisson flow is the flow $\phi_t$ on $P/G$
induced by $\Phi_t$. The orbit $x_e=Gp_e\in P/G$ is an equilibrium
point of $\phi_t$. Moreover the relative equilibrium $p_e$ is
$G$-stable if and only if $x_e$ is Lyapunov stable in the usual
sense, and is leafwise stable in the sense of
Definition~\ref{REStabDef} if and only if $x_e$ is leafwise stable in
the sense of~\S\ref{PoiStabSect}. We can therefore expect to lift the
theory developed in~\S\ref{PoiStabSect} to obtain both topological and
derivative tests for the $G$-stability of relative
equilibria.

The topological and derivative tests for the $G$-stability of
Hamiltonian relative equilibria are contained in Sections
\ref{HamStabTestSect} and \ref{T2EMCM}, respectively.
Section \ref{transverse} introduces
\emph{transverse Poisson structures} as tools for describing local
leafspace topology and uses this to discuss some special cases.
At the end of Section \ref{T2EMCM} we discuss applications
to relative equilibria of systems that are invariant under actions of
the Euclidean groups $\SE(2)$ and $\SE(3)$ and an
example of Libermann and Marle \cite{LibermannPMarleCM-1987.1}.

\subsection[$G$-Stability by Topological Methods]{$\mathbf G$-Stability by Topological Methods}
\label{HamStabTestSect}

In this section we obtain topological stability criteria for relative
equilibria of $G$-invariant Hamiltonians on $P$
by applying the results of Section \ref{PoiTopSect}
for Poisson equilibria to the symmetry reduced flow on the
Poisson orbit manifold $P/G$.

The main result of the section is the following theorem.
It is a generalisation to non-compact free group actions of
Montaldi's stability theorem~\cite{MontaldiJA-1997.1} for compact groups.
An extension to general proper actions is given in~\cite{PRW2}.

\begin{theorem}
\label{HamTopThm}

Let $H$ be a $G$-invariant Hamiltonian on $P$ with a
relative equilibrium at $p_e$. Let $S$ be a slice transverse to
the orbit $G p_e$ at $p_e$. Then $p_e$ is
$G$-stable if there exists an open neighbourhood $U_S$ of $p_e$ in
$S$ and an open neighbourhood $U_{\mu_e}$ of $\mu_e$ in $\mfk g^*$
such that
\begin{equation}\label{HamTomCdn}
H^{-1}\bigl(H(p_e)\bigr)\cap J^{-1}\bigl(T_2^{U_{\mu_e}}(\mu_e)\bigr)
\cap U_S=\{p_e\}.
\end{equation}
This remains true if $H$ is replaced by any conserved
quantity with values in any Hausdorff space.
\end{theorem}

\noindent
The proof of Theorem \ref{HamTopThm} uses the slice $S$ as a local
model for the orbit space $P/G$. The projection from
$S$ to $P/G$ induces a Poisson structure on $S$ that is isomorphic to that
on the corresponding open neighbourhood of $Gp_e$ in $P/G$.
The following result shows that the $T_2$-set
of the Poisson structure on $S$ at $p_e$ is just the pullback of the
$T_2$-set of the Poisson structure on $\g^*$ at $\mu_e$.
\begin{lemma}
\label{T2S}
Let $S$ denote a slice to $Gp_e$ at $p_e$.
\begin{enumerate}
\item
The symplectic leaves of the induced Poisson structure on
$S$ are the connected components of the intersections
$J^{-1}({\cal O}) \cap S$,
where $\cal O$ is coadjoint orbit in $\g^*$.
\item
There exist arbitrarily small neighbourhoods $U_S$ of $p_e$ in $S$
and $U_{\mu_e}$ of $\mu_e$ in $\g^* $ such that
\begin{equation}
\label{T2slice}
T_2^{U_S}(p_e) = J^{-1}\bigl(T_2^{U_{\mu_e}}(\mu_e)\bigr)\cap U_S.
\end{equation}
\end{enumerate}
\end{lemma}

\bpr
The first part follows from the fact that the symplectic leaves of
$P/G$ are the connected components of $J^{-1}({\cal O})/G$.

For the second part we note that since $J$ is a submersion there
exist arbitrarily small neighbourhoods $U_{p_e}$ of $p_e$ in $P$
and $U_{\mu_e}$ of $\mu_e$ in $\g^*$ such that if $M$ is a connected
set in $U_{\mu_e}$ then $J^{-1}(M) \cap U_{p_e}$ is also connected.
If $U_S = U_{p_e}\cap S$ then this also implies that
$J^{-1}(M) \cap U_S$ is connected, since this last intersection
is homeomorphic to the projection of $J^{-1}(M) \cap U_{p_e}$
to $P/G$.
It follows that the momentum map induces a homeomorphism between
the leaf space of $U_S$ and that of $U_{\mu_e}$, which in turn implies
(\ref{T2slice}).
\epr

\noindent
{\bf Proof of Theorem \ref{HamTopThm}.} \quad
We may assume that the neighbourhoods $U_S$ and
$U_{\mu_e}$ in (\ref{HamTomCdn}) also satisfy the conclusions of
Lemma \ref{T2S}. The condition (\ref{HamTomCdn}) therefore implies that the
restriction $H\vert_S$ isolates $p_e$ in $T_2^{U_S}(p_e)$ and so
Theorem \ref{PoiStabTh1} implies that $p_e$ is Lyapunov stable for the flow
on $S$ generated by $H\vert_S$ and the induced Poisson structure on
$S$. Identifying $S$ with a neighbourhood of $Gp_e$ in $P/G$ and the
flow on $S$ with the quotient of the flow on $P$ generated by $H$ gives
the result.

The statement that $H$ can be replaced by any conserved quantity with values
in any Hausdorff space follows immediately from the analogous statement
in Theorem \ref{PoiStabTh1}.
\hfill $\ink$

\begin{remark}\label{leafstab}
It follows from Lemma \ref{T2S} that
the leaf space of $P/G$ is Hausdorff at $Gp_e$
if and only if $\g^*/G$ is Hausdorff at $G\mu_e$.
If this holds then the condition
(\ref{HamTomCdn}) reduces to
\begin{equation}\label{HamLeafCdn}
H^{-1}\bigl(H(p_e)\bigr)\cap J^{-1}(\mu_e)\cap U_S=\{p_e\},
\end{equation}
which means that $H$ isolates $G_{\mu_e}p_e$ in its momentum level set.
Sufficient conditions on $\mu_e$ for
$\g^*/G$ to be Hausdorff at $G\mu_e$
are given in Proposition \ref{prop:reg+split}.
In general condition (\ref{HamLeafCdn}) always implies that $p_e$ is
leafwise stable.
\end{remark}

\subsection{Transverse Poisson Structures and Tame Generators}
\label{transverse}

Weinstein's local splitting theorem for a Poisson manifold
$X$ \cite[Theorem 2.1]{WeinsteinA-1983.1}
 states that any point $x \in X$ has
a neighbourhood $U$ that is isomorphic as a Poisson manifold
to the product of a neighbourhood of $x$ in the leaf $L(x)$
and a transverse Poisson space $\Sigma$.
The isomorphism class of this transverse Poisson structure
does not depend on either the point $x$ in $L(x)$ or the choice of
transverse section $\Sigma$. In Section \ref{transverseg*} we
give a detailed description of the transverse Poisson space
at a point $\mu \in \g^*$ and use this to give sufficient conditions
for $\g^*/G$ to be locally Hausdorff at $G_{\mu}$ and for
a generator $\xi \in \g_{\mu}$ to be tame. In
Section \ref{transverseP/G} we show that the transverse Poisson
structure at $Gp \in P/G$ is isomorphic to that at $\mu = J(p)$ and use
this to describe the tame generators on $P/G$ at $Gp$.

\subsubsection{Coadjoint Actions}
\label{transverseg*}

In this section we will define an explicit choice of a
transverse Poisson space to the coadjoint orbit
$G\mu$ at $\mu$ in $\g^*$
and describe some of its properties.

Let $G_\mu$ denote the isotropy subgroup of $\mu\in\g^*$ and $\g_\mu$
the isotropy subalgebra. Let $\n_\mu$ denote a complement to $\g_\mu$
in $\g$, so that $\g=\g_\mu\oplus\n_\mu$ and $\g^*=\g_\mu^\ro
\oplus\n_\mu^\ro$, where the superscript~$^{\rm o}$ denotes an
annihilator. The decomposition of $\g$ induces an isomorphism
between $\n_\mu^\ro$ and $\g_\mu^*$. An easy calculation shows that
$\mfk g\mu=\mfk g_\mu^\ro$, so the affine subspace $\mu+\n_\mu^\ro$ is
transverse to the coadjoint orbit $G\mu$ at $\mu$. It then follows
from~\cite[Proposition 4.1]{WeinsteinA-1983.1} that there is an
induced Poisson
structure on a neighbourhood $\Sigma$ of $\mu$ in $\mu+\n_\mu^\ro$.
This transverse Poisson structure can be described explicitly as follows
\cite{CR01,RWL-1999}. Let $\pi_{\g_\mu^\ro}$ denote
the projection from $\g^*$ to $\g_\mu^\ro$ with kernel $\n_\mu^\ro$.
Since $\Sigma$ is affine, its tangent spaces are all equal to
$\mfk n_\mu^\ro$,
which is identified with $\g_\mu^*$. Thus the derivatives of functions
on $\Sigma$ can be regarded as taking values in $\g_\mu$.

\begin{proposition}\label{SigmaProp}{\rm\cite{CR01,RWL-1999}}
There exists a neighbourhood $V\subseteq\n_\mu^\ro$ of $0$ such that:
\begin{enumerate}
\item For every $\nu \in V$ and every $\xi\in\g_\mu$ the equation
\begin{equation}
\pi_{\g_\mu^\ro}\left(\ad^*_{\xi+\eta}(\mu+\nu)\right)=0
\end{equation}
has a unique solution $\eta=\eta_{\mu+\nu}(\xi)\in\n_\mu$ and the map
$j_\mu(\nu):\g_\mu\to\g$ defined by
$j_\mu(\nu)\xi=\xi+\eta_{\mu+\nu}(\xi)$ is linear, depends smoothly on
$\nu$ and satisfies $j_\mu(0)\xi = \xi$.
\item
A Poisson structure is defined on $\Sigma = \mu + V$ by the bracket
\begin{equation*}
\{f,g\}(\mu+\nu)=-\bigl\langle\mu+\nu,
[j_\mu(\nu)d_\nu f(\mu+\nu), j_\mu(\nu)d_\nu g(\mu+\nu)]\bigr\rangle
\end{equation*}
where $f,g$ are smooth functions on $\Sigma$ and $[\cdot,\cdot]$ is the
Lie bracket on $\g$.
\item
The Poisson structure on a neighbourhood of $\mu$ in $\g^*$ is isomorphic
to the product of this Poisson structure on $\Sigma$ and the
Kostant-Kirillov-Souriau symplectic structure on a neighbourhood of
$\mu$ in $G\mu$.
\end{enumerate}
\end{proposition}

The symplectic leaves of the Poisson structure of $\Sigma$ are the
intersections of $\Sigma$ with the coadjoint orbits on $\g^*$. Next
we characterise these using some notation and results from~\cite{W01}.
\begin{definition}\label{ZmunuDef}
For each $\nu$ in a neighbourhood $ V\subseteq\n_\mu^\ro$ of $0$
define $Z_{\mu,\nu}$ to be the connected component containing the identity
of
\begin{equation*}
\tilde{Z}_{\mu,\nu}\equiv\set{g\in G}{\Ad^*_{g}(\mu+\nu)\in\mu + V},
\end{equation*}
and define $Z_{\mu,V}=\bigcup_{\nu \in V} Z_{\mu,\nu}$.
Clearly $Z_{\mu,0} = G_\mu^\ro$ and $Z_{\mu,\nu}\supseteq
G_{\mu+\nu}^\ro$.
\end{definition}

\begin{proposition}\label{ZmunuProp}\mbox{}
\begin{enumerate}
\item
For each $\nu$ sufficiently close to $0$ there exists a neighbourhood
$W$ of the identity in $G$ such that $Z_{\mu,\nu} \cap W$ is a manifold
of dimension $\dim G_\mu$. Its tangent space at the identity is
$T_{\rm id}Z_{\mu,\nu} = j_\mu(\nu)\g_\mu$.
\item
For $V$ sufficiently small the symplectic leaf through $\mu + \nu$
of the transverse Poisson structure on $\Sigma = \mu + V$ is equal to
 $Z_{\mu,\nu}(\mu+\nu) \cap \Sigma$.
\end{enumerate}
\end{proposition}
The second statement says essentially that the symplectic leaves of
$\Sigma$ are the `orbits' of $Z_{\mu,V}$. However in general
$Z_{\mu,\nu}$ and $Z_{\mu,V}$ are not groups.
\bpr
The first statement was proved in~\cite{W01} and follows from an application
of the implicit function theorem. The second statement is an immediate
consequence of Definition~\ref{ZmunuDef} and the fact that the symplectic
leaves of $\Sigma$ are its intersections with the coadjoint orbits of
$G$.\epr

The set of generators (see
Definition~\ref{Ptame}) at $\mu$ for the
Poisson structure on $\g^*$ is $(\mfk g\mu)^\ro=\mfk g_{\mu}$. Since any element of $T_{\mu}^*\Sigma$ is a
generator for the Poisson structure on $\Sigma$, and
$T_{\mu}^*\Sigma=\n_{\mu}^*= (\mfk g\mu)^\ro=\mfk g_{\mu}$, the generators
for the two Poisson structures coincide.

By Weinstein's local splitting theorem, for an appropriate
neighbourhood $U_{\mu}$ of $\mu\in\g^*$, the set $T_2^{U_{\mu}}(\mu)$
is the product of
$T_2^\Sigma(\mu)$ with a neighbourhood of $\mu$ in the leaf $G\mu$.
Thus, a generator $\xi\in\g_{\mu}$ is {\it tame} as a generator on
$\g^*$ if and only if it is tame as a generator on $\Sigma$,
ie if and only if the equivalent conditions
$T_{\mu}\bigl(T_2^\Sigma(\mu)\bigr)\subseteq\onm{ann}_{\g_\mu^*}\xi$
and $T_{\mu}\bigl(T_2^{U_{\mu}}(\mu)\bigr)\subseteq\onm{ann}_{\g^*}\xi$
hold. The annihilator of $\xi$ in $\g_\mu^*$ is denoted by
$\onm{ann}_{\g_\mu^*}\xi$, and other annihilators in a similar way.
Note that the property of being tame does not depend on the choice of
the transverse section $\Sigma$.

In the next section we will need a stronger form of tameness.
\begin{definition}
\label{verytame}
A generator $\xi\in T_\mu^*\Sigma\cong\g_\mu$ is {\it very tame}
if $T_2^\Sigma(\mu)\subseteq\mu+\onm{ann}\xi$.
\end{definition}
Clearly very tame generators are always tame,
but in general the converse does not hold, see Remark \ref{nonsplitex}.

There are two important special cases of transverse Poisson structures.
\begin{definition}
We say that $\mu\in\g^*$ is:
\begin{enumerate}
\item
{\it Regular} if $\dim\g_\nu=\dim \g_\mu$ for every $\nu$ in a
neighbourhood of $\mu$.
\item
{\it Split} if there exists a $G_\mu^\ro$-invariant complement
$\n_\mu$ to $\g_\mu$ in $\g$.
\end{enumerate}
\end{definition}
Here $G^\ro$ denotes the identity component of $G$.
Note that $\mu$ is regular if and only if it is a regular
point for the Lie-Poisson structure on $\g^*$ in the sense of
Section \ref{PoiDerSect},
since the symplectic leaves are the coadjoint orbits and the
definition above implies that these have constant dimension
near regular points.

It is clear that if there exists a $G$-invariant inner product
on $\g$, as in the case
of compact or Abelian groups, then every $\mu$ is split,
the $T_2$-sets are trivial and every generator is very tame.
Some other results are given in the following proposition.

\begin{proposition}\label{prop:reg+split}\mbox{}
\begin{enumerate}
\item
If $\mu$ is regular then $\g_\mu$ is Abelian, $Z_{\mu,\nu} =
G^\ro_{\mu+\nu}$ for $\nu$ small, the transverse Poisson structure on
$\Sigma$ is trivial, $\g^*/G$ is Hausdorff at $G\mu$ and
every generator is very tame. The set of
regular points in $\g^*$ is open and dense.
\item
Let $\mu$ be split. Then $Z_{\mu,\nu} = G_\mu^\ro$ for every $\nu$, the
transverse Poisson structure on $\Sigma$ is isomorphic to the standard
Lie-Poisson structure on $\g_\mu^*$ and its symplectic leaves are the
coadjoint orbits of $G_\mu^\ro$ on $\g_\mu^*$. Moreover a generator
$\xi\in\g_\mu$ is tame if and only if it is very tame, and a sufficient
condition
for $\xi\in\g_\mu$ to be tame is that its adjoint orbit $G_\mu^\ro\xi$ is
bounded.
Finally $\g^*/G$ is Hausdorff at $G\mu$ if there is a
$G_\mu^\ro$-invariant
inner product on $\g_\mu$ and in particular if $G_\mu$ is compact.
\end{enumerate}
\end{proposition}

\bpr The results for regular $\mu$ are due essentially to
Lie~\cite{L1890} and Duflo and Vergne~\cite{DV69}. For split $\mu$
the results in the first sentence of Part~\emph{2} follow from
Proposition~\ref{SigmaProp}. See also~\cite{WeinsteinA-1983.1, CR01}
and references therein. For the remaining statements of Part~\emph{2} note that
since $\mu$ is split we can assume that $\mu=0$,
$G_\mu=G$ and $\Sigma$ is an open neighbourhood of $0$ in $\mfk
g^*$. The linearity of the action of $G$ on $\g^*$ implies that if
$\nu\in T_2^\Sigma(0)$ then $\lambda\nu\in T_2^\Sigma(0)$ for all
$\lambda \in \R$. It follows that $\nu\in
T_0\bigl(T_2^\Sigma(0)\bigr)$ and so $T_2^\Sigma(0)\subseteq
T_0\bigl(T_2^\Sigma(0)\bigr)$. Hence, if $\xi$ annihilates
$T_0\bigl(T_2^\Sigma(0)\bigr)$ then it annihilates $T_2^\Sigma(0)$,
and so the sets of tame and very tame generators coincide.

If $\nu \in T_2^{\Sigma}(0)$ then there are sequences $\nu_n\in\Sigma$
and $g_n\in G_\mu^\ro$, such that $\nu_n\to \nu$ and $g_n \nu_n\to 0$
as $n\to\infty$. We need to show that $\langle\nu,\xi\rangle=0$ if
$G^\ro\xi$ is bounded. But $\langle\nu,\xi\rangle=\lim_{n\to\infty}
\langle\nu_n,\xi\rangle=\lim_{n\to\infty}\langle g_n\nu_n,g_n
\xi\rangle=0$ since $g_n\nu_n\to0$, $g_n\xi\in G^\ro\xi$ and
$G^\ro\xi$ is bounded.\epr

\begin{remark} \label{nonsplitex}
In general regular momenta need not be split and their $Z_{\mu,\nu}$
need not be equal to $G_\mu^\ro$. Examples include nonzero points on
the nilpotent cone of $sl(2,\R)^*$ (see Example \ref{sl(2)} and
Example 4.4 of \cite{RWL-1999}). If $\mu$ is non-split and
non-regular then there may be tame generators $\xi\in \g_\mu$ which
are not very tame. In Example \ref{nosmoothing}, where $\mu$ is a
subregular nilpotent element of $sl(3,\R)^*$, the fact that
$T_\mu\bigl(T_2^{U_\Sigma}(\mu)\bigr)$ is one dimensional while
$T_2^{U_\Sigma}(\mu)$ is two dimensional implies that there exist
generators $\xi$ which annihilate
$T_\mu\bigl(T_2^{U_\Sigma}(\mu)\bigr)$ but not $T_2^{U_\Sigma}(\mu)$
itself.
\end{remark}

\subsubsection[The Orbit Space $P/G$]{The Orbit Space $\mathbf{P/G}$}
\label{transverseP/G}

We now give a description of how the transverse
Poisson structure at $\mu \in \g^*$ is related to local Poisson structure
of $P/G$ near $Gp$.

As in Section \ref{HamStabTestSect} we identify a neighbourhood of
$Gp$ in $P/G$ with a slice $S$ transverse to the orbit $Gp$ in $P$. By
Lemma \ref{T2S} the symplectic leaf $L^S(p)$ of the Poisson structure
on $S$ is the connected component of $J^{-1}(G\mu)\cap S$ containing
$p$.  Since $J$ is a submersion the image $\Sigma=J(S)$ is a
submanifold of $\g^*$ through $\mu$ that is transverse to $G\mu$ and
so has an induced transverse Poisson structure. The fact that $J$ is
a Poisson map \cite{MarsdenJERatiuTS-1994.1} implies that that the
transverse Poisson space at $p\in S$ is isomorphic to that at
$\mu\in\mfk g^*$. It follows from Weinstein's local splitting theorem
that a sufficiently small neighbourhood of $p$ in $S$ is isomorphic as
a Poisson manifold to a neighbourhood of $(p,\mu)$ in the product
$L^S(p)\times\Sigma$.

Recall that $p_e$ is a relative equilibrium with generator $\xi_e$
if and only if it is a critical point of $H_{\xi_e}$.
The derivative tests which we give in Section \ref{T2EMCM}
state that if the generator
is tame (or very tame) and the Hessian
$d^2H_{\xi_e}(p_e)$ is definite when restricted to certain subspaces
of the normal space $N$ to $Gp_e$ at $p_e$, then $p_e$ is $G$-stable.
The {\it Witt decomposition} \cite{BaLe97} of the normal space to
any group orbit $Gp$ splits it into two components
$N = N_0 \oplus N_1$, where $N_1$ is the {\it symplectic normal space},
ie a maximal subspace of $N$ for which the restriction of the symplectic
form on $T_{p}P$ is nondegenerate. More explicitly $N_1$
is a complement to $T_{p}(G_{\mu}p)$ in $\ker dJ(p)$.
For free group actions it is isomorphic to the tangent
space at $Gp$ to the leaf $L(Gp)$ through $Gp$ in $P/G$, ie the
Marsden-Weinstein reduced phase space \cite{MW74}.

The derivative $dJ(p_e)$ maps a
complement $N_0$ to $N_1$ in $N$ isomorphically to $T_{\mu_e}\Sigma$,
the normal space to $G\mu$ at $\mu$ and hence
identifies it with $\n_{\mu}^\ro \cong \g_{\mu}^*$.
Using this identification it can be
endowed with the transverse Poisson structure at $\mu$ described in
Section \ref{transverseg*}. The normal space $N$ can then be given
the Poisson structure obtained by taking the product of this structure
on $N_0$ with the symplectic structure on $N_1$.

The following result is a version of Weinstein's
local splitting of $P/G$, or equivalently a slice $S$,
due to Guillemin, Sternberg and Marle.

\begin{theorem} {\rm \cite{GS84,M85}}
\label{MGSTh}
Let $S$ be a slice in $P$ transverse to $Gp$ at $p$ such that
$J(S) \subset \Sigma = \mu + \n_{\mu}^\ro$. Then
there exists a Poisson isomorphism from 
$S$ to an open neighbourhood of the
origin in $N = N_0 \oplus N_1$ with its product Poisson structure
such that the restriction of the momentum map $J$ to
$S$ is given by $J_S(\nu,w) = \mu + \nu$ where
$\nu \in N_0 \cong \n_{\mu}^\ro$ and $w \in N_1$.
\end{theorem}

It follows from the local decomposition of $P/G$ given in
Theorem \ref{MGSTh} that the set of generators at $Gp_e$ can
also be identified with $T_{\mu}^*\Sigma$, and hence with $\mfk g_{\mu}$,
and that a generator is (very) tame for the Poisson structure on $P/G$ if
and only if under this identification it is (very) tame for the
transverse Poisson structure on $\Sigma \cong \g_{\mu}^*$.

\subsection[$T_2$-Energy-Momentum-Casimir Method]{$\mathbf{T_2}$-Energy-Momentum-Casimir Method}
\label{T2EMCM}

We now turn to the task of lifting the derivative tests
of Section \ref{PoiDerSect} from $P/G$ to the symplectic phase space $P$.
Our main results are Theorem \ref{HamHessStabTh1} and
Corollary \ref{Cor:Proceedings} which, when combined with the results of
Section \ref{AStability}, generalise all previous published results for free
group actions. Their extensions to general proper actions will be given in
\cite{PRW2}.

Before giving our main result we present a simple corollary of
Theorem \ref{HamTopThm}. This will be generalised by Corollary
\ref{Cor:Proceedings}.

\begin{corollary}\label{HamHessCor}
Let $p_e$ be a relative equilibrium of $H$ with generator $\xi_e$.
Suppose that $\g^*/G$ is Hausdorff at $\mu_e = J(p_e)$
and that the Hessian $d^2 H_{\xi_e}(p_e)$ is (positive or negative)
definite when restricted to any symplectic normal space at $p_e$.
Then $p_e$ is $G$-stable.
\end{corollary}

\noindent
If $\g^*/G$ is not Hausdorff at $\mu_e$
then the
definiteness of $d^2 H_{\xi_e}(p_e)$
on a symplectic normal space
only implies that $p_e$ is leafwise stable.

By Proposition \ref{prop:reg+split} the space
$\g^*/G$ is Hausdorff at $\mu_e$
if either $\mu_e$ is regular or
$\mu_e$ is split and there
exists a $G_{\mu_e}$-invariant inner product on $\g_{\mu_e}^*$.
The $G$-stability result for $\mu_e$ regular is
Theorem~8.17 of Chapter~IV of~\cite{LibermannPMarleCM-1987.1}.
The condition that $\mu_e$ is split and there
exists a $G_{\mu_e}$-invariant inner product on $\g_{\mu_e}^*$
is implied by the existence of a $G_{\mu_e}$-invariant inner
product on the whole of $\g^*$, which is always true for compact
groups. Under this stronger hypothesis
\cite{LeSi98,OrRa99} show that definiteness of $d^2 H_{\xi_e}(p_e)$
on a symplectic normal space implies that $p_e$ is actually
$G_{\mu_e}$-stable, generalising a result of \cite{P92} for compact
groups $G$. We recover this result by combining Corollary
\ref{HamHessCor} with Corollary \ref{compact_A} in
Section \ref{AStability}.

\bpr If $S$ is a slice to $Gp_e$ at $p_e$ the symplectic normal space
can be identified with $T_{p_e}\bigl(J^{-1}(\mu_e)\cap S\bigr)$ and
the definite Hessian, together with the fact that $p_e$ is a critical
point of $H_{\xi_e}$, implies that $p_e$ is isolated in
$H_{\xi_e}^{-1}\bigl(H_{\xi_e}(p_e)\bigr) \cap J^{-1}(\mu_e)\cap S$.
Since $J_{\xi_e}$ is constant on $J^{-1}(\mu_e)$ this implies that
$p_e$ is also isolated in $H^{-1}\bigl(H(p_e)\bigr) \cap
J^{-1}(\mu_e)\cap S$. If $\g^*/G$ is Hausdorff at $\mu_e$ then this
is equivalent to condition (\ref{HamTomCdn}) of Theorem
\ref{HamTopThm} and so $p_e$ is also $G$-stable. \epr

In Theorem \ref{HamHessStabTh1} below we generalise the
`energy-momentum method' of Corollary \ref{HamHessCor}
to obtain criteria for the $G$-stability of Hamiltonian
relative equilibria for general non-compact symmetry groups and
non-regular momentum values.
Before formulating the theorem we make some preliminary remarks
on smoothings and Casimirs.
A~\emph{local Casimir} is a continuous
function $C:U_{\mu_e}\to \R$
on a neighbourhood $U_{\mu_e}$ of $\mu_e$
which is constant on
the connected components of the coadjoint orbits intersected with
$U_{\mu_e}$.

\begin{remark}\label{B_iECMRem}\mbox{}

\begin{enumerate}
\item
It follows from Lemma \ref{T2S} that for a sufficiently small slice
$S$ to $Gp_e$ at $p_e$ we have
$T_2^S(p_e) = J^{-1}(T_2^\Sigma(\mu_e)) \cap S$
where $\Sigma = J(S)$.
Since $J$ is a submersion at $p_e$ this means that
$\{B_i\}_{i=1}^n$ is a smoothing of $T_2^S(p_e)$
if and only if $\{J(B_i)\}_{i=1}^n$ is a smoothing of
$T_2^\Sigma(\mu_e)$.
\item
Every local Casimir $C:U_{\mu_e}\to \R$ with $\Sigma \subset U_{\mu_e}$
is constant on $T_2^\Sigma(\mu_e)$.
By Lemma \ref{T2S} the function $C \circ J$
is a Casimir on every slice $S$ satisfying $J(S) \subseteq U_{\mu_e}$.
Consequently, as in Section \ref{PoiStabSect}, if
$\{B_i\}_{i=1}^n$ is a smoothing of $T_2^S(p_e)$
the function $C\circ J|_{B_i }$ has a critical point at $p_e$ for each
$i$.
\item
If $\{B_i\}_{i=1}^n$ is a smoothing of $T_2^S(p_e)$
a relative equilibrium $p_e$ with generator $\xi_e \in \g_{\mu_e}$ is
a critical point of $H|_{B_i }$ for each $i$ if and only if $\xi_e$ is tame.
To see this, note that $dH(p_e)=\langle dJ(p_e),\xi_e\rangle$
annihilates
$T_{p_e} T_2^S(p_e) = dJ(p_e)^{-1}\bigl(T_{\mu_e}T_2^\Sigma(\mu_e)\bigr)$
if and only if $\xi_e$ annihilates
$T_{\mu_e}\bigl(T_2^\Sigma(\mu_e)\bigr)$.
\end{enumerate}
\end{remark}

\noindent
The next theorem provides two Hessian tests for the $G$-stability of
relative equilibria with tame generators. The first uses the second
derivatives of the restrictions of $H$ to the manifolds $B_i$ in a
smoothing of $ T^S_2(p_e) = J^{-1}(T_2^\Sigma(\mu_e)) \cap S$.
In general it is
not possible to calculate the second derivative of $H|_S$ and then
restrict it to the tangent spaces $T_{p_e}B_i$ because $H$ need not
have a critical point at $p_e$.
If $J^{-1}\bigl(T^\Sigma_2(\mu_e)\bigr)=J^{-1}(\mu_e)$
then this can be circumvented by replacing the
second derivative of $H$ by that of $H_{\xi_e}=H-J_{\xi_e}$, which
does have a critical point at $p_e$. This is valid because $J_{\xi_e}$ is
constant on this level set of $J$.
We used this idea in the proof of Corollary \ref{HamHessCor}.

If the set $J^{-1}\bigl(T^\Sigma_2(\mu_e)\bigr)$ is larger
than the level set of $J$ then it
may not be possible to replace $H$ by $H_{\xi_e}$. However if $\xi_e$
is {\it very tame}, as defined in Definition \ref{verytame},
then $J_{\xi_e}$ is constant on the whole of
$J^{-1}\bigl(T^\Sigma_2(\mu_e)\bigr)$ and it is possible to proceed as in
the Hausdorff case. This gives the second test. Both tests can be
strengthened by the addition of Casimirs, as stated in the theorem.

\begin{theorem} {\rm ($T_2$-Energy-Momentum-Casimir Method)}
\label{HamHessStabTh1}
Let $p_e$ be a relative equilibrium of $H$ with generator $\xi_e$.
Let $S$ be a slice to $Gp_e$ at $p_e$ and $\{B_i\}_{i=1}^n$
a smoothing of $T_2^S(p_e)$.
\begin{enumerate}
\item
If $\xi_e$ is a tame generator then $p_e$ is $G$-stable
if for each $1\le i\le n$ there is a smooth local Casimir $C_i$
such that the Hessian $d^2\bigl((H+C_i\circ J)|B_i\bigr)(p_e)$ is
definite on $T_{p_e} B_i$.
\item
 If $\xi_e$ is a very tame generator for the transverse Poisson
structure on $\Sigma = J(S)$ then $p_e$ is $G$-stable if for
each $1\le i\le n$ there is a smooth local Casimir $C_i$
with a critical point at $\mu_e$ and
a Hessian $d^2\bigl(H_{\xi_e}+C_i\circ J\bigr)(p_e)$ that is definite
when restricted to $T_{p_e} B_i$.
\end{enumerate}
\end{theorem}
\bpr

\noindent{\em 1.}\quad
By Remark \ref{B_iECMRem}.2 the functions $C_i \circ J$ are Casimirs on $S$
and so {\em 1.} follows from
Theorem \ref{PoiHessStabTh}.

\smallskip\noindent{\em 2.}\quad
If $\xi_e$ is very tame then for all
$p\in T_2^S(p_e) = J^{-1}\bigl(T_2^\Sigma(\mu_e)\bigr)$
we have:
\begin{equation*}
J_{\xi_e}(p)=\langle J(p),\xi_e\rangle
=\langle J(p)-\mu_e,\xi_e\rangle+\langle\mu_e,\xi_e\rangle
=\langle\mu_e,\xi_e\rangle
\end{equation*}
since $J(p)\in T_2^\Sigma(\mu_e)\subseteq\mu_e+\onm{ann}\xi_e$.
Let $f_i = H|_S + C_i \circ J|_S$, $\tilde{f}_i = H_{\xi_e}|_S + C_i \circ
J|_S$,
$f = (f_1,\ldots, f_n)$ and $\tilde{f} = ( \tilde{f}_1,\ldots, \tilde{f}_n)$.
Then on $J^{-1}\bigl(T_2^\Sigma(\mu_e)\bigr)$
the restrictions of $f$ and $\tilde{f}$ differ by a constant vector.

If $C_i$ has a critical point at $\mu_e$ then $\tilde{f}_i$ has a
critical point at $p_e$ and
\begin{equation*}
d^2\tilde{f}_i(p_e)|_{T_{p_e}B_i}
 =d^2( \tilde{f}_i|{B_i})(p_e)
\end{equation*}
It follows that if the restriction of the Hessian of $\tilde{f}_i$ to
$T_{p_e}B_i$ is definite then $p_e$ is an isolated point in
$(\tilde{f}_i|B_i)^{-1}\bigl(\tilde{f}_i(p_e)\bigr)$ for each $i$ and
so an isolated point in $\tilde{f}^{-1}\bigl(\tilde{f}(p_e)\bigr)\cap
J^{-1}\bigl(T_2^\Sigma(\mu_e)\bigr)$. Since on
$J^{-1}\bigl(T_2^\Sigma(\mu_e)\bigr)$ the functions $f$ and
$\tilde{f}$ differ only by a constant this implies that $p_e$ is an
isolated point of $f^{-1}\bigl(f(p_e)\bigr)\cap
J^{-1}\bigl(T_2^\Sigma(\mu_e)\bigr)$. The proof now follows from
Theorem \ref{HamTopThm}.\epr 
We will extend the $T_2$-Energy-Momentum-Casimir Method to general
proper actions in~\cite{PRW2}.

\begin{remark}\label{ECMRem}\mbox{}

\begin{enumerate}
\item
By Remark \ref{B_iECMRem}.1
inclusion of Casimirs in the theorem is
unnecessary if the smoothing $\{J(B_i)\}_{i=1}^{n}$ of $T_2^\Sigma(\mu_e)$
satisfies conditions similar to those stated at the end of~\S\ref{PoiDerSect}.
In particular Casimirs are unnecessary if $T_2^\Sigma(\mu_e)$ is a manifold,
and the inclusion of $C_i$ is unnecessary for any $i$ such that
(\ref{casUnnecEq}) holds for the Poisson space $\Sigma$.
\item
If $T_2^S(p_e)$ is a one dimensional
submanifold of $S$ then it is not difficult to formulate
a derivative test for the $G$-stability
of $p_e$ analogous to that of Corollary~\ref{PoiStabTh2}
for Poisson equilibria.
\item
Theorem \ref{HamHessStabTh1} can be extended to weak smoothings,
as in Remark \ref{weakSmoothRem}.
\end{enumerate}
\end{remark}

We announced the following corollary of Theorem \ref{HamHessStabTh1}
in \cite[Theorem 2]{PRWProc}. Although the result is not always
optimal, it gives a generalisation of the energy-momentum method of
\cite{LeSi98,P92,OrRa99} to non-compact symmetry groups that, as will
be seen in Corollaries \ref{SE2cor} and \ref{SE3cor}, is optimal for
the Euclidean symmetry groups that are most likely to arise in
applications. Recall from Section \ref{transverseP/G} that the normal
space $N$ to $Gp_e$ at $p_e$ decomposes as $N_0 \oplus N_1$ where
$N_1$ is a symplectic normal space and $N_0 \cong \g_\mu^*$.

\begin{corollary}
\label{Cor:Proceedings}{\rm (Simple Energy-Momentum Method)}
Let ${\mfk t}_{\mu_e}\subseteq\g_{\mu_e}$ be the space of very tame
generators at $\mu_e$
and ${\mfk w}^*_{\mu_e}=\onm{ann}_{\g_{\mu_e}^*}{\mfk t}_{\mu_e}$.
Then the relative equilibrium $p_e$ with momentum $\mu_e$
is $G$-stable if its generator $\xi_e$ is very tame and
$d^2 H_{\xi_e}(p_e)|_{{\mfk w}_{\mu_e}^* \oplus N_1}$ is definite.
\end{corollary}
We call the subspace
${\mfk w}^*_{\mu_e}= \onm{ann}_{\g_{\mu_e}^*}{\mfk t}_{\mu_e}$ of
$\g_{\mu_e}^*$ the space of {\em wild momenta}.

\bpr
Choose the slice $S$ at $p_e$ as in Theorem \ref{MGSTh}
and identify it with the normal space $N = N_0 \oplus N_1$.
By definition very tame generators annihilate the set $T_2^\Sigma(\mu_e)$
and so the space $B ={\mfk w}^*_{\mu_e} \oplus N_1 \subseteq N$ is a weak
smoothing of $T_2^S(p_e)$.
The corollary therefore follows from Remark \ref{ECMRem}.3.
\epr

\noindent
For the case of semidirect products of compact groups and vector
spaces Leonard and Marsden \cite{LeonardNEMarsdenJE-1997.1} also identified an
intermediate set between $N_1$ and $T_pP$ on which definiteness of
$d^2H_{\xi_e}(p_e)$ implies stability. Corollary \ref{Cor:Proceedings}
generalises their results.

In \cite{PRWProc} we restricted to the case where $\mu_e$ is split
and {\it defined} generators $\xi \in \g_{\mu_e}$ to be
tame if $G_{\mu_e}^\ro\xi$ is bounded.
By Proposition \ref{prop:reg+split}.2 this definition of tame generators
is stronger than the notion of very tame generators used in this paper.

We now show how the energy-momentum method of this section applies to
relative equilibria of Euclidean invariant Hamiltonian systems. We
treat the cases where the symmetry group is the special Euclidean
group $G=\SE(2)=\SO(2)\ltimes\R^2$ of rotations and translations of
the plane and the special Euclidean group $G=\SE(3)=\SO(3)\ltimes
\R^3$ of three-space. Relative equilibria of Euclidean invariant
Hamiltonian systems have been studied for several systems, including
the dynamics of underwater vehicles
\cite{LeonardNE-1997.1,LeonardNEMarsdenJE-1997.1} and systems of point
vortices \cite{Newton}, and the statics of elastic rods \cite{Mielke}.

In the case of $\mtl{SE}(2)$ the coadjoint action is
\begin{equation*}
(R,a)\cdot(\mu^r,\mu^a)=(\mu^r+R\mbb J\mu^a\cdot a,R\mu^a)
\end{equation*}
where we have represented elements of $\mtl{SO}(2)$ as $2\times2$
matrices and $\mbb
J=\left[\begin{array}{cc}0&1\\-1&0\end{array}\right]$. From this one
sees that $(\mu^r,\mu^a)\mapsto|\mu^a|$ is a Casimir, and that the
coadjoint orbits are the cylinders about the $\mu^r$ axis together
with the points on the $\mu^r$ axis. Consequently
$\mu\in\mtl{se}(2)^*$ is regular if $\mu^a\ne0$ in which case
$G_\mu=\set{(1,t\mu_e^a)}{t\in\mbb R}$, and $\mu\in\mtl{se}(2)^*$ is
non-regular if $\mu^a=0$ in which case $G_\mu=\mtl{SE}(2)$. Any
regular momentum $(\mu^r,\mu^a)$ is split since the transverse section
$\set{(\mu^r,t\mu^a)}{t\in\mbb R}$ is $G_{\mu}$ invariant, while the
non-regular momenta are trivially split since they are at zero
dimensional coadjoint orbits. Any two neighbourhoods of two
non-regular momenta both meet sufficiently narrow cylinders about the
$\mu^r$-axis, so if $\mu$ is non-regular then $T_2(\mu)$ is the
$\mu^r$-axis. Given a non-regular momentum $(\mu^r,\mu^a)$, a
generator $(\xi^r,\xi^a)$ is tame if it annihilates the tangent space
to this $T_2$-set or equivalently if and only if $\xi^r=0$, and every
tame generator is very tame. Casimirs are not required for an
application of Theorem~\ref{HamHessStabTh1} in this case since the
$T_2$-sets are manifolds. 
The $T_2$-set at the nonregular  momenta is equal to the annihilator of
the (very) tame generators and so  Theorem~\ref{HamHessStabTh1} is
equivalent to Corollary~\ref{Cor:Proceedings} here. 
Summarising all
this gives the following corollary.
\begin{corollary}\label{SE2cor}
Assume $G=\mtl{SE}(2)$. Let $p_e$ be a relative equilibrium of $H$
with generator $\xi_e=(\xi^r_e,\xi^a_e)$ and momentum
$\mu_e=(\mu^r_e,\mu^a_e)$. 
\begin{enumerate}
\item If $\mu_e^a\ne 0$ then $\mu_e$ is regular,
$G_{\mu_e}=\set{(1,t\mu_e^a)}{t\in\mbb R}$, and $p_e$ is $G$-stable if
$d^2H_{\xi_e}(p_e)$ is definite on 
any complement to $\mfk g_{\mu_e}p_e$ in $\onm{ker} dJ(p_e)$ (ie definite
on any symplectic normal space).
\item If $\mu_e^a=0$ then $\mu_e$ is non-regular,
$G_{\mu_e}=\mtl{SE}(2)$, $\xi_e$ is tame if and only if $\xi_e^r=0$,
and $p_e$ is $G$-stable if $\xi_e$ is tame and $d^2H_{\xi_e}(p_e)$ is
definite on any complement to $\mfk gp_e$ in
$dJ^{-1}(p_e)(\mtl{so}(2)^*)$.
\end{enumerate}
\end{corollary}

The coadjoint action of $\mtl{SE}(3)$ is explicitly computed in
\cite{MarsdenJERatiuTS-1994.1} and the analysis proceeds in a similar
way to the $\mtl{SE}(2)$ case. Again the 
$T_2$-sets are subspaces of $\mtl{se}(3)^*$, so that every
tame generator is also very tame, even though the momenta which
are nonzero and non-regular are not split. The result is the following
corollary.
\begin{corollary}\label{SE3cor}
Assume $G=\mtl{SE}(3)$. Let $p_e$ be a relative equilibrium of $H$
with generator $\xi_e=(\xi^r_e,\xi^a_e)$ and momentum
$\mu_e=(\mu^r_e,\mu^a_e)$. 
\begin{enumerate}
\item If $\mu_e^a\ne 0$ then $\mu_e$ is regular, 
$G_{\mu_e}=
\set{(R,t\mu_e^a)}{R\mu_e^a=\mu_e^a,t\in\mbb R}
\cong\mtl{SO}(2)\times\mbb R$,
and $p_e$ is $G$-stable if
$d^2H_{\xi_e}(p_e)$ is definite on 
any complement to $\mfk g_{\mu_e}p_e$ in $\onm{ker}dJ(p_e)$ (ie definite
on any symplectic normal space).
\item If $\mu_e^a=0$ then $\mu_e$ is non-regular, 
\begin{equation*}
G_{\mu_e}=\begin{cases}
\set{(R,a)}{R\mu_e^r=\mu_e^r,a\in\mbb R^3}\cong\mtl{SO}(2)\ltimes\mbb R^3&\mu_e^r\ne 0,\\
\mtl{SE}(3)&\mu_e^r=0,
\end{cases}\end{equation*}
$\xi_e$ is tame if and only if $\xi_e^r=0$,
and $p_e$ is $G$-stable if $\xi_e$ is tame and $d^2H_{\xi_e}(p_e)$ is
definite on any complement to $\mfk gp_e$ in
$dJ^{-1}(p_e)(\mtl{so}(3)^*)$.
\end{enumerate}
\end{corollary}

For both the $\mtl{SE}(2)$ and $\mtl{SE}(3)$ symmetry groups our
theory indicates that, in the case of zero translational momentum,
only for purely translating relative equilibria of $\mtl{SE}(2)$
invariant systems is $G$-stability accessible by energy-momentum
methods, since only then is the generator tame. For spinning relative
equilibria which have zero translational momentum $G$-stability may be
accessible in low dimensional phase space by a blow up argument
coupled with KAM theory, as we will show in \cite{PRW3}.

\begin{example}\label{RsdR}

This example, which is due to Libermann, Marle, and Krishnaprasad, is
Exercise~(15.10) on page~274 of \cite{LibermannPMarleCM-1987.1} and is
worked out in detail in the appendix
of~\cite{KrishnaprasadPS-1989}. It was presented as an example of a
relative equilibrium which is leafwise stable, but not $G$-stable.
Here we show how it relates to the theory of this paper.

The group $G$ is a semidirect product $G =\mbb R\ltimes\mbb R$
with multiplication
\begin{equation*}
(a,b)(a^\prime,b^\prime)=(a+a^\prime,b+e^ab^\prime)
\end{equation*}
The phase space $P$ is $T^*G=\mbb R^4=\{(q^1,q^2,p_1,p_2)\}$
with its standard symplectic form $dq^1\wedge dp_1 + dq^2\wedge dp_2$.
The action of $G$ on $P$ is the cotangent lift of the left translation
action of $G$ on itself:
\begin{equation*}
(a,b)(q^1,q^2,p_1,p_2)=(a+q^1,b+e^aq^2,p_1,e^{-a}p_2).
\end{equation*}
The Lie algebra dual is $\mfk g^*$ is $\mbb R^2=\{(\nu_1,\nu_2)\}$,
the coadjoint action is
\begin{equation*}
\onm{Ad}^*_{(a,b)^{-1}}\left(\begin{array}{c}\nu_1\\\nu_2\end{array}\right)
=
\left(\begin{array}{cc}1&e^{-a}b\\0&e^{-a}\end{array}\right)
\left(\begin{array}{c}\nu_1\\\nu_2\end{array}\right),
\end{equation*}
and the $\Ad^*$-equivariant momentum mapping is
\begin{equation*}
J(q^1,q^2,p_1,p_2)=
\left(\begin{array}{c}p_1+p_2q^2\\p_2\end{array}\right).\end{equation*}
As is evident from the formula for the coadjoint action, there are
zero dimensional coadjoint orbits along the $\nu_1$ axis, and the half
planes $\nu_2>0$ and $\nu_2<0$ are two dimensional coadjoint orbits
of regular momenta. The $T_2$-sets for the non-regular momenta on
the $\nu_2$ axis are $\mfk g^*$, and consequently a relative
equilibrium at such a momentum is tame if and only its generator is
zero (ie if and only if it is an equilibrium).

The Hamiltonian considered in
\cite{LibermannPMarleCM-1987.1,KrishnaprasadPS-1989}
is $H = p_2e^{q^1}$, for which the flow on $P$ is
\begin{equation*}
q^1(t)=q^1(0),\quad q^2(t)=q^2(0)+te^{q^1(0)},\quad
p_1(t)=p_1(0)-t p_2(0)e^{q^1(0)},\quad p_2(t)=p_2(0).
\end{equation*}
The relative equilibria are given by solving $dH=dJ_{(\dot a,\dot b)}$,
which directly leads to the equations
\begin{equation*}
p_2e^{q^1}=0,\quad 0=\dot ap_2,\quad 0=\dot a,\quad e^{q^1}=\dot aq^2+\dot b,
\end{equation*}
with solution
\begin{equation*}
p_2=0,\quad\dot a=0,\quad\dot b=e^{q^1}.
\end{equation*}
The momenta of these relative equilibria all have $\nu_2=0$ and hence
are all non-regular, and the generators all have $\dot b\ne 0$, and
hence are all wild, so we do not expect that any of the relative
equilibria are $G$-stable. On the other hand the momentum level sets
of the non-regular momenta are all group orbits so the
Marsden-Weinstein reduced spaces corresponding to those momenta are
all points, and consequently all the relative equilibria are leafwise
stable.

To explicitly see that the relative equilibria are not $G$-stable, 
it is best to
examine the flow on the Poisson space $P/G$. The standard
Lie-Poisson reduction to $\mfk g^*$ is realized by the
quotient map given by the left trivialisation
$\nu=(\nu_1,\nu_2)=g^{-1}(p_1,p_2)=(p_1,p_2e^{q^1})$.
The flow descends to
\begin{equation*}
\nu_1(t)=\nu_1(0)-t\nu_2(0),\quad \nu_2(t)=\nu_2(0)
\end{equation*}
and the relative equilibria descend to the equilibria $\nu_2=0$. These
equilibria are not stable by direct inspection of the reduced flow:
perturbations into $\nu_2\ne 0$ of such equilibria are carried far
from their origin by translation parallel to the $\nu_1$ axis.
\end{example}

\section{Beyond $\mathbf G$-Stability}
\label{AStability}
In this section we show that $G$-stable relative equilibria typically
satisfy a stronger stability property. To motivate this, suppose
$p(t)$ is an integral curve of $X_H$ starting at $p(0)\approx p_e$. If
$p_e$ is $G$-stable then there is a curve $g(t)$ in $G$ such that
$g(t)^{-1}p(t)\approx p_e$. By continuity and by equivariance and
conservation of the momentum $J$ we have:
\begin{equation*}
\mu_e=J(p_e)\approx J\bigl(g(t)^{-1}p(t)\bigr)
=\onm{Ad}^*_{g(t)}J\bigl(p(t)\bigr)
=\onm{Ad}^*_{g(t)}J\bigl(p(0)\bigr)\approx\onm{Ad}^*_{g(t)}\mu_e.
\end{equation*}
In other words \emph{conservation of momentum ought to imply that
trajectories which start close to a $G$-stable relative equilibrium
$p_e$ should remain close to the $G_{\mu_e}$ orbit through $p_e$,
and not just the $G$ orbit.} Results of this type on
$G_{\mu_e}$-stability, all of which assume Hessian conditions and
compactness-related invariant inner products or norms on $\mfk g$,
first appeared in~\cite{P92}, with extensions to non-free actions
in~\cite{LeSi98,Matsui-2001.1,OrRa99}.
That $G$-stable relative equilibria are not generally
$G_{\mu_e}$-stable in the non-compact case was noticed in~\cite{P92} and then
in~\cite{LeonardNEMarsdenJE-1997.1}.
In this paper, assuming only the $G$-stability of a relative
equilibrium, we obtain results on stability which lie between
$G_{\mu_e}$-stability and $G$-stability. The results are proved for
free actions, but extend to general proper actions \cite{PRW2}.

In Section \ref{AStabilityThm} we give the main general theorem.
This is specialised to split momenta in Section \ref{Asplit}, and then
applied to Euclidean invariant Hamiltonian systems in Section
\ref{AEuclidean}.

\subsection[$A$-Stability]{$\mathbf A$-Stability}
\label{AStabilityThm}

We first define a very general stability property for equivariant flows.
\begin{definition}
\label{Adefn} Let $P$ be a topological space with an action of a group
$G$ and a $G$-equivariant flow $\Phi_t\colon P\to P$. For any subset $A
\subseteq G$, a point $p_e\in P$ is said to be $A$-{\em stable} with
respect to the flow $\Phi_t$ if for any open neighbourhood $
U\subseteq P$ of $p_e$ there exists an open neighbourhood $ U^\prime\subseteq
P$ of $p_e$ such that if $p\in U^\prime$ then $\phi_t(p)\in AU$ for all
$t\in\mbb R$.
\end{definition}
If $A\subseteq B\subseteq G$ then if $p_e$ is $A$-stable it
is also $B$-stable. It is easy to show that a point $p_e$ is $A$-stable if it is $AW$-stable for any neighbourhood 
$W$ of $1\in G$.

Let $\Sigma\subseteq\mu_e+\n_{\mu_e}^\ro$ be a transverse section
through $G\mu_e$ at $\mu_e = J(p_e)$.
For any small neighbourhood $V$ of $0$ in
$\n_{\mu_e}^\ro$ and any point $\nu\in V$ let $Z_{\mu_e,\nu}$ and
$Z_{\mu_e,V}$ be the subsets of $G$ as in
Definition~\ref{ZmunuDef}. For any neighbourhood $W$ of $1$ in $G$
define
\begin{equation*}
A_{V,W}(\mu_e)=\bigcup_{g\in W}gZ_{\mu_e,V}g^{-1}.
\end{equation*}
Since $Z_{\mu_e,0}=G^\ro_{\mu_e}$ we have $G^\ro_{\mu_e} \subseteq
A_{V,W}(\mu_e)\subseteq G$ for any pair of neighbourhoods $V$ and 
$W$.
\begin{theorem}
\label{AStabThm}
If $p_e$ is a $G$-stable relative equilibrium
then $p_e$ is an $A_{V,W}(\mu_e)$-stable equilibrium for every
pair of neighbourhoods $ V $ of $0$ in $\n_{\mu_e}^\ro$ and
$ W$ of $1$ in $G$.
\end{theorem}

\bpr
Let $S$ be a slice at $p_e$ which is mapped by $J$ to the transverse
section $\Sigma \subseteq\mu_e+\n_{\mu_e}^\ro$ at $\mu_e$.
By Theorem \ref{MGSTh} the slice $S$ can be identified
with the product of open neighbourhoods of $\{0\}$ in $N_0 $ and $N_1$,
and in these coordinates the restriction of $J$ to $S$ is given
by $J_S(\nu,w)=\mu_e+\nu$ where $\nu\in N_0=\n_{\mu_e}^\ro$. By
equivariance a $G$-invariant neighbourhood of $Gp_e$ can be parametrised
by $G \times S$ and in these coordinates the momentum map is given
by $J(g,\nu,w)=\Ad^*_{g^{-1}}(\mu_e+\nu)$.

Let $\Phi_t$ denote the $G$-equivariant flow on $P$ and $\phi_t$ the
induced flow on $P/G$. As in Section \ref{HamStabTestSect}, we identify $S$
with an open neighbourhood of $Gp_e$ in $P/G$ and denote the flow
on $S$ by $\phi_t$. Let $U$ be a sufficiently small
neighbourhood of $p_e$.
We then have $W\times(\mu_e+U_0) \times U_1 \subseteq U$ where $W$ is a
neighbourhood of $1\in G$, and $U_0$ and $U_1$ are
neighbourhoods of $\{0\}$ in $N_0$ and $N_1$, respectively.

Since $Gp_e$ is Lyapunov stable for the flow $\phi_t$ on $P/G$ there
exists a neighbourhood $U_S^\prime$ of $p_e$ in $U_S=(\mu_e+U_0)\times U_1$
such that $\phi_t(U_S^\prime)\subseteq U_S$ for all $t$.
We may choose $U_S^\prime$ to be a product $(\mu_e+U_0^\prime)\times
U_1^\prime$ where $U_0^\prime$ and $U_1^\prime$ are again neighbourhoods of
$\{0\}$ in $N_0$ and $N_1$, respectively. Since $S$ is a
slice, for $g_0\in G$, $\nu_0\in U_0^\prime$ and $w_0\in U_1^\prime$ we have
\begin{equation*}
\Phi_t(g_0,\nu_0,w_0)\equiv\bigl(g(t),\nu(t),w(t)\bigr)=
g(t)\bigl(1,\nu(t),w(t)\bigr)
\equiv g(t)\phi_t(\nu_0,w_0).
\end{equation*}
The momentum map $J$ is preserved by the flow $\Phi_t$ and so
$\Ad^*_{g(t)^{-1}}\bigl(\mu_e+\nu(t)\bigr)=\Ad^*_{g_0^{-1}}\bigl(\mu_e+\nu_0\bigr)$
for all $t$, and hence
$\Ad^*_{g_0^{-1}g(t)}(\mu_e+\nu_0)=\mu_e+\nu(t)\in\mu_e+U_0$. If
the neighbourhood $U_0\subseteq V$ is chosen small enough this implies that
$g(t) \in g_0Z_{\mu_e,\nu_0}$ for all $t$ (see Definition~\ref{ZmunuDef}).
It follows that
\begin{equation*}
\Phi_t(g_0,w_0,\nu_0)=\gamma(t)\big(g_0,w(t),\nu(t)\bigr)=\gamma(t)\bigl(g_0,
 \phi_t(w_0,\nu_0)\bigr)
\end{equation*}
where $\gamma(t)=g(t)g_0^{-1}\in g_0Z_{\mu_e,\nu_0}g_0^{-1}$.
Hence $\Phi_t(W\times U_S)\subseteq A_{W,V}(\mu_e)U$, as required.
\epr

\subsection{Split Momenta}
\label{Asplit}

Theorem \ref{AStabThm} can be improved when $\mu_e$ is split.
In this case Proposition~\ref{prop:reg+split} says
$Z_{\mu_e,V} = G^\ro_{\mu_e}$ and so
$A_{V,W}(\mu_e)= G^{\ro,W}_{\mu_e}
\equiv\bigcup_{g \in W} gG^\ro_{\mu_e}g^{-1}$.
Let $G^0_{\mu_e} = L_{\mu_e}K_{\mu_e}$ where $L_{\mu_e}$
is a submanifold of $G^\ro_{\mu_e}$ and
$K_{\mu_e}$ a subgroup of $G^0_{\mu_e}$ for which there exists a
$K_{\mu_e}$-invariant inner product on $\g^*$.
Such a splitting exists for any connected Lie group with
$K_{\mu_e}$ a maximal compact subgroup of $G^\ro_{\mu_e}$
\cite[Theorem 3.1]{Hochschild}.
However in some cases it may be possible to take $K_{\mu_e}$ larger
than this.
For any neighbourhood $W$ of $1$ in $G$ define
$L_{\mu_e}^W = \bigcup_{g \in W} gL_{\mu_e}g^{-1}$.
\begin{corollary}\label{splitmu_A}
If $p_e$ is a $G$-stable relative equilibrium and $\mu_e = J(p_e)$
is split then $p_e$ is an $L_{\mu_e}^WK_{\mu_e}$- stable equilibrium
for every neighbourhood $W$ of $1$ in $G$.
\end{corollary}
\bpr
Let $\bar K_{\mu_e}$ denote the quotient of $K_{\mu_e}$ which acts
effectively on $\g^*$. The kernel $T$ of the homomorphism $K_{\mu_e} \to
\bar K_{\mu_e}$ acts trivially on $\g^*$, and hence also on $\g$. It
follows that $T$ must lie in the centre of $G^\ro$.

For any neighbourhood $ W $ of $1$ in $G$ define
$K_{\mu_e}^W=\bigcup_{g\in W} gK_{\mu_e}g^{-1}.$
 For any
neighbourhood $ U$ of $p_e$ in $P$ we claim that there exist
neighbourhoods $W$ of $1$ in $ G$ and $U^\prime$ of $p_e$ in $ P$
such that $K_{\mu_e}^W U^\prime\subseteq K_{\mu_e}U$. Indeed consider
the continuous map $\tau\colon G^\ro\times\bar K_\mu \times P\to P$
defined by $\tau(g,kT,p)=k^{-1}gkg^{-1}p$. Note that this is
well defined because $T$ lies in the centre of $G^\ro$. So given $U$ and
$kT$ there exist open neighbourhoods $W_k$ of $1$ in $G^\ro$, $V_k$
of $kT$ in $\bar K_{\mu_e}$ and $U_k$ of $p$ in $P$ such that
$\tau(W_k,V_k,U_k) \subseteq U$. Since $\bar K_{\mu_e} $ is compact
there exists a finite subcover $\{V_{k_i}\}_{i=1}^a$ of $\bar K_{\mu_e} $.
Put $W = \bigcap_{i=1}^aW_{k_i}$ and $U^\prime=\bigcap_{i=1}^aU_{k_i}$
to prove the claim.

For such a choice of $U'$ and $W$ we have
\[
G^{\ro,W}_{\mu_e} U' \ \subseteq\ L_{\mu_e}^W K_{\mu_e}^W U'
\ \subseteq \ L_{\mu_e}^W K_{\mu_e} U.
\]
It follows that if $p_e$ is $G^{\ro,W}_{\mu_e}$-stable then
it is also $L_{\mu_e}^W K_{\mu_e}$-stable.
\epr
The following corollary generalises the results of \cite{P92, LeSi98, OrRa99}
by requiring only that $p_e$ is $G$-stable and not that the
restriction of $d^2H_{\xi_e}(p_e)$ to the symplectic normal space $N_1$
is definite.
\begin{corollary}
\label{compact_A}
If $p_e$ is a $G$-stable relative equilibrium and there exists a
$G^\ro_{\mu_e}$-invariant inner product on $\g^*$ then
$p_e$ is $G^\ro_{\mu_e}$-stable.
\end{corollary}

We end this section by applying Corollary \ref{splitmu_A}
to the case of a semidirect product $G=K\ltimes \V$ with
$K$ a compact Lie group acting
linearly on a finite dimensional real vector space $\V$.
Assume that $\mu_e$ is split and $G_{\mu_e}$ has the form
$G_{\mu_e} = K_{\mu_e} \ltimes \V_{\mu_e}$ where $K_{\mu_e}$ is a
closed subgroup of $K$ and $\V_{\mu_e}$ is a
$K_{\mu_e}$ invariant subspace of $\V$. Then $G_{\mu_e} = \V_{\mu_e}
K_{\mu_e}$.
Let $W$ denote an open neighbourhood of $(1,0)$ in $K \ltimes \V$
of the form $W_K \times W_\V$ where $W_K$ is a neighbourhood of $1$
in $K$ and $W_\V$ is a neighbourhood of $0$ in $\V$. Then
a calculation shows $\V_{\mu_e}^W = W_K^{-1}\V_{\mu_e} \subseteq \V$,
a generalised `cone' in $\V$. Set theoretically $A = W_K^{-1}\V_{\mu_e} K_{\mu_e}$ is equal to
$K_{\mu_e} \times W_K^{-1}\V_{\mu_e} \subseteq K \times \V$, so we
obtain the following result.

\begin{corollary}
\label{semidirect_A}
Let $p_e$ be a $K\ltimes\V$-stable relative equilibrium for which $\mu_e$ is
split and $G_{\mu_e} = K_{\mu_e} \ltimes \V_{\mu_e}$. Then $p_e$ is a
$K_{\mu_e} \times W_K^{-1}\V_{\mu_e}$-stable relative equilibrium for any open
neighbourhood $W_K$ of $1$ in $K$.
\end{corollary}

\subsection{Euclidean Invariant Hamiltonian Systems}
\label{AEuclidean}

Finally, we apply the $A$-stability results of this section to relative
equilibria of Euclidean invariant Hamiltonian systems and briefly
discuss their implications for the stability of rigid bodies in fluids
\cite{LeonardNE-1997.1,LeonardNEMarsdenJE-1997.1}. As in Corollaries
\ref{SE2cor} and \ref{SE3cor} we treat the cases where the symmetry
group is the special Euclidean group $G=\SE(2)$ of the plane or the
special Euclidean group $G=\SE(3) $ of three-space. We assume that
the relative equilibria are $G$-stable and compute subsets $A\subset
G$ for which the relative equilibria are $A$-stable.

First consider the symmetry group $G=\SE(2)$. For a regular momentum
value $\mu_e$ (case 1 of Corollary \ref{SE2cor}) let $C$ be any open
cone containing $\mu_e^a$ and $A=\{1\}\times C$. Then, for some open
neighbourhood $W_{\SO(2)}$ of $1\in\mtl{SO}(2)$,
$C\subseteq\bigcup_{R\in W_{\SO(2)}} R\mbb R\mu_e^a$ and so $p_e$ is
$A$-stable by Corollary \ref{semidirect_A}. The cone can be made
arbitrarily `thin' but smaller cones will require initial conditions
closer to $p_e$. For a non-regular momentum value (case 2 of Corollary
\ref{SE2cor}) we have $G_{\mu_e}=G$ so trivially $A = G$ since $A$
must contain $G_{\mu_e}$.

Now let $G=\SE(3)$. As in the case of $\mtl{SE}(2)$ symmetry,
for a regular momentum value $\mu_e$ (case 1 of
Corollary \ref{SE3cor}) we can apply Corollary \ref{semidirect_A}
and conclude that $p_e$ is $A$-stable for
$A=\mtl{SO}(2)\times C$ where $C$ is any open cone containing
$\mu_e^a$ and $\SO(2)$ consists of rotations about $\mu_e^a$. For zero
momentum $A=\mtl{SE}(3)$ trivially. 
Case 2 of Corollary \ref{SE3cor} with $\mu_e^r\ne0$ is more interesting because $\mu_e$
is not split and so Corollary \ref{semidirect_A} does not apply and a
direct analysis of the sets $A_{V,W}(\mu_e)=\bigcup_{g\in
W}gZ_{\mu_e,V}g^{-1}$ used in Theorem \ref{AStabThm} is required. 
The coadjoint
orbit through $(\mu_e^r,0)$ is two sphere
$\set{(\mu^r,0)}{|\mu^r|=|\mu_e^r|}$ and so a transverse section to
the orbit can be taken as
\begin{equation*}
\Sigma=\set{\mu_e+(\nu^r,\nu^a)}{\nu^r=t\mu^r_e,t\in\mbb R, \nu^a\in\mbb R^3}
\cong\mbb R\times\mbb R^3.
\end{equation*}
Let  $\theta_{R,\mu_e^r}$  denote the  angle  between $R\mu_e^r$  and
$\mu_e^r$.  We claim that, if $\epsilon_0>0$ and $\epsilon_1>0$, and
\begin{equation*}
A_{\epsilon_0,\epsilon_1}\equiv
\set{(R,a)\in\mtl{SE}(3)}{|\sin\theta_{R,\mu_e^r}
|<\epsilon_1|a|+\epsilon_0},
\end{equation*}
then there are neighbourhoods $V\subseteq\Sigma$ of $\mu_e$ and $W$ of
$1\in\mtl{SE}(3)$ such that $A_{V,W}\subseteq
A_{\epsilon_0,\epsilon_1}$. Therefore by Theorem \ref{AStabThm} an $\mtl{SE}(3)$
stable relative equilibrium with a momentum as in Case~2 of Corollary \ref{SE3cor} with $\mu_e^r\ne0$ is
$A_{\epsilon_0,\epsilon_1}$-stable for any $\epsilon_0>0$ and
$\epsilon_1>0$.
To prove the claim,
using the $\mtl{SE}(3)$ coadjoint action \cite{MarsdenJERatiuTS-1994.1}
one computes that $Z_{\mu_e,\nu}$ is the set of $(R,a)\in\mtl{SE}(3)$
satisfying
\begin{equation}\label{e1}
\mu_e^r\times\bigl(R(\mu_e^r+\nu^r)+a\times R\nu^a\bigr)=0.
\end{equation}
In this equation we think that $\nu^r=t\mu_e^r$ where $t$ and $\nu^a$ are
small. If $a$ is small then $\mu_e^r\times R\mu_e^r\approx0$, so $R$
is nearly a rotation about $\mu_e^r$, but the equation admits
solutions for $R\in\mtl{SO}(3)$ arbitrary as long as $a$ is large
enough. We rewrite (\ref{e1}) as
\begin{equation*}
(1+t)\mu_e^r\times R\mu_e^r=-\mu_e^r\times(a\times R\nu^a)
\end{equation*}
whereupon
\begin{equation}\label{e5}
|\sin\theta_{R,\mu_e^r}|=\frac{|\mu_e^r\times R\mu_e^r|}{|\mu_e^r|^2}\le
\frac{|\nu^a|}{|\mu_e^r|}|a|.
\end{equation}
The action of conjugation is 
\begin{equation*}
(R^\prime,a^\prime)\equiv(\tilde R,\tilde a)(R,a)(\tilde R,\tilde a)^{-1}
=(\tilde R R\tilde R^{-1},\tilde Ra-\tilde R R\tilde R^{-1}\tilde a+\tilde a).
\end{equation*}
By continuity, for arbitrary
$\epsilon_0$ there is a neighbourhood $W$ of $1\in\mtl{SE}(3)$ such
that if $(\tilde R,\tilde a)\in W$ then
\begin{equation}\label{e6}
|\sin\theta_{R^\prime,\mu_e^r}|=
|\sin\theta_{\tilde RR\tilde R^{-1},\mu_e^r}|
<|\sin\theta_{R,\mu_e^r}|+\frac{\epsilon_0}2.
\end{equation}
The required $R$ uniformity of this estimate follows from compactness
of $\mtl{SO}(3)$. Combining (\ref{e5}) and (\ref{e6}), and choosing
$V$ so that $|\nu^a|<|\mu_e^r|\epsilon_1$ gives
\begin{equation*}
|\sin\theta_{R^\prime,\mu_e^r}|
<\frac{|\nu^a|}{|\mu_e^r|}|a|
+\frac{\epsilon_0}2<\epsilon_1|a|+\frac{\epsilon_0}2.
\end{equation*}
Also, by shrinking $W$ there is the uniform (in $a$ and
$R$) estimate
\begin{equation*}
|a^\prime|=|\tilde Ra-\tilde R R\tilde R^{-1}\tilde a+\tilde a|
\ge|a|-|\tilde R R\tilde R^{-1}\tilde a-\tilde a|
\ge|a|-\frac{\epsilon_0}{2\epsilon_1}
\end{equation*}
from which follows
\begin{equation*}
|\sin\theta_{R^\prime,\mu_e^r}|<\epsilon_1\left(
|a^\prime|+\frac{\epsilon_0}{2\epsilon_1}\right)+\frac{\epsilon_0}2
=\epsilon_1|a^\prime|+\epsilon_0,
\end{equation*}
as required.

This $A_{\epsilon_0,\epsilon_1}$-stability implies orientation
stability only when translation is a priori confined. Notice that for
the set $A$ given by Corollary \ref{splitmu_A} for split $\mu_e$ the
projection of $A$ to the compact part of $G$ is the same as that of
$G_{\mu_e}$. Equation \ref{e1} shows that in the non-split case these
projections can be very different.

The results described in this section were inspired by the stability
analysis of Leonard and Marsden \cite{LeonardNEMarsdenJE-1997.1} of
the Kirchhoff model for the dynamics of a rigid body in a fluid. For
$G$ a semidirect product of a compact group and a vector space they
prove $\Gamma$-stability results for certain groups $\Gamma$ between
$G_{\mu_e}$ and $G$. For regular momentum values of $\SE(2)$-invariant
systems their group $\Gamma$ is $\R^2$, while for regular momentum
values of $\SE(3)$-invariant systems it is $\SO(2)\ltimes\R^3$. However
numerical integration of an example with $\SE(3)$ symmetry suggests
that the drift in the translational direction only occurs within a
cone, in accord with the results of this section.

\section*{Acknowledgments}
This work was partially supported by an EPSRC Visiting Fellowship
(GR/L57074) and an NSERC individual research grant for GWP, an EPSRC
Research Grant (GR/K99893), a Scheme Four grant from the London
Mathematical Society, a European Community Marie Curie Fellowship
(HPMF-CT-2000-00542) for CW, and by European Community funding for the
Research Training Network `MASIE' (HPRN-CT-2000-00113). The authors
thank the University of Warwick Mathematics Institute for its
hospitality during several visits when parts of the paper were
written.

\footnotesize\frenchspacing
\bibliographystyle{plain}

\end{document}